\documentclass[12pt]{amsart}

\usepackage{amssymb,latexsym,amscd}

\numberwithin{equation}{section}

\newcommand{\doublesubscript}[3]{
\displaystyle\mathop{\displaystyle #1_{#2}}_{#3} }

\newtheorem{theorem}{Theorem}[section]
\newtheorem{proposition}[theorem]{Proposition}
\newtheorem{conjecture}[theorem]{Conjecture}
\newtheorem{corollary}[theorem]{Corollary}

\newtheorem{lemma}[theorem]{Lemma}

\newtheorem{maintheorem}[theorem]{Main Theorem}
\theoremstyle{definition}

\newtheorem{remark}[theorem]{Remark}

\newtheorem{definition}[theorem]{Definition}

\def\ZZ{\mathbb{Z}}

\def\GG{\mathbb{G}}

\def\kk{\Bbbk}

\def\gg{\mathfrak{g}}

\def\hh{\mathfrak{h}}

\def\GG{\mathcal G}

\def\NN{\mathcal{N}}

\def\FF{\mathcal{F}}

\advance\oddsidemargin by-0.5in
\advance\evensidemargin by-0.5in
\advance\textwidth by 1in

\addtocounter{section}{-1}

\begin{document}

\title{Lie algebras and Lie groups  over noncommutative rings}

\author{Arkady Berenstein}
\address{Department of Mathematics, University of Oregon,
Eugene, OR 97403} \email{arkadiy@math.uoregon.edu}

\author{Vladimir Retakh}
\address{\noindent Department of Mathematics, Rutgers University, Piscataway, NJ 08854}
\email{vretakh@math.rutgers.edu}

\thanks{The authors were supported in part
by the NSF grant DMS \#0501103 (A.B.), and  by the NSA grant H98230-06-1-0028 (V.R.).
}

\makeatletter
\renewcommand{\@evenhead}{\tiny \thepage \hfill  A.~BERENSTEIN and  V.~RETAKH \hfill}

\renewcommand{\@oddhead}{\tiny \hfill Lie algebras and Lie groups over  noncommutative rings
 \hfill \thepage}
\makeatother

\date{May 27, 2007}

\begin{abstract}
The aim of this paper is to introduce and study  Lie algebras and
Lie groups over noncommutative rings. For any Lie algebra $\gg$
sitting inside an associative algebra $A$ and any associative
algebra $\FF$ we introduce and study the algebra $(\gg,A)(\FF)$,
which is the Lie subalgebra  of $\FF \otimes A$ generated by $\FF
\otimes \gg$. In many examples $A$ is the universal enveloping
algebra of $\gg$. Our description of the algebra $(\gg,A)(\FF)$
has a striking resemblance to the commutator expansions of $\FF$
used by M. Kapranov in his approach to noncommutative geometry.
To each algebra $(\gg, A)(\FF)$ we associate a ``noncommutative
algebraic'' group which naturally acts on $(\gg,A)(\FF)$ by
conjugations  and conclude the paper with some examples of such
groups.
\end{abstract}

\maketitle

\tableofcontents

%\vspace{-.25in}

\pagebreak[3]

\section{Introduction}

The aim of this paper is to introduce and study algebraic groups
and Lie algebras over noncommutative rings. 

Our approach is motivated by the following considerations.
A naive definition of a Lie algebra as a bimodule over a
noncommutative associative algebra $\FF$ (over a field $\kk$) does not bring
any interesting examples beyond Lie algebra
$gl_n(\FF)$. Even the special Lie
%$gl_n(\FF)=M_n(\FF)=\FF\otimes gl_n(\kk)$. Even the special Lie
algebra $sl_n(\FF)=[gl_n(\FF), gl_n(\FF)]$
%$sl_n(\FF)=[gl_n(\FF), gl_n(\FF)]$ (which is the Lie
%subalgebra of all matrices in $gl_n(\FF)$ whose traces belong to the commutator $[\FF,\FF]$)
is not an $\FF$-bimodule. Similarly, the special  linear group
$SL_n(\FF)$ is not defined by equations but rather by congruences
given by the Dieudonne determinant (see \cite{A}). This is why
the ``straightforward'' approach to classical groups over rings
started by J.~Dieudonne in \cite{D} and continued by
O.~T.~O'Meara and others (see \cite {HO}) does not lead to
new algebraic groups. Also, unlike in the commutative case, these
methods do not employ rich structural theory of Lie algebras.

As a starting point, we observe that the Lie algebra $sl_n(\FF)$ where $\FF$ is an
associative algebra over a field $\kk$ (of characteristic $0$) 
%usually defined as the kernel of the trace map $M_n(\FF)\rightarrow \FF/[\FF,\FF]$, 
is the Lie subalgebra of $M_n(\FF)=\FF \otimes M_n(\kk)$ generated by $\FF\otimes sl_n(\kk)$ 
(all tensor products in the paper are taken over $\kk$ unless
specified otherwise).
This motivates us to consider, for any Lie subalgebra $\gg$ of an
associative algebra $A$, the Lie subalgebra
\begin{equation}
\label{eq:pivlie}
(\gg,A)(\FF)\subset \FF\otimes A
\end{equation}
generated by $\FF\otimes \gg$.

If $\FF$ is commutative, then $\FF\otimes \gg$ is already a Lie
algebra, and  $(\gg,A)(\FF)=\FF\otimes \gg$. However, if $\FF$ is noncommutative, 
this equality does not hold. Our first main result (Theorem \ref{th:perfect}) is a
formula expressing $(\gg,A)(\FF)$ in terms of powers
$\gg^n=Span\{g_1g_2\dots g_n: g_1,\dots ,g_n\in \gg\}$ in $A$ for
all {\it perfect pairs} $(\gg, A)$
in the sense of Definition \ref{def:perfect pair}. The class of perfect pairs is large enough -- 
it includes all  semisimple and Kac-Moody Lie algebras $\gg$.

More precisely, Theorem \ref{th:perfect} states that for any perfect pair $(\gg,A)$  and any 
associative algebra $\FF$ we have
\begin{equation}
\label{eq:mainintro1}
(\gg,A)(\FF)=\FF\otimes \gg+\sum _{n\geq 1}I_n\otimes [\gg, \gg^{n+1}]
+[\FF, I_{n-1}]\otimes \gg^{n+1} \ ,
\end{equation}
where $I_0\supset I_1\subset I_2\supset \cdots $ is the descending filtration of two-sided ideals of $\FF$ defined inductively as
$I_0=\FF$, $I_{n+1}=\FF[\FF, I_n]+[\FF, I_n]$. It is remarkable that this filtration emerged in works by M. Kapranov in \cite {K} and
then by M.~Kontsevich and A.~Rosenberg in \cite{KR} as an important tool
in noncommutative geometry. 

We also get a compact formula (Theorem \ref{th:perfect sl_2}) for $(sl_2(\kk),A)(\FF)$ which, hopefully, has physical implications.

For unital $\FF$ we define the {\it ``noncommutative'' current group} or, in short, 
the $\NN$-current group  $G_{\gg,A}(\FF)$ to be the set of all
invertible elements $X\in \FF\otimes A$  such that $X\cdot (\gg, A)(\FF)\cdot X^{-1}=(\gg, A)(\FF)$.
This is our generalization of $GL_n(\FF)$. In fact, if $\gg=sl_n(\kk)\subset M_n(\kk)=A$ then $G_{\gg,A}(\FF)=GL_n(\kk)$.
 
However, for other compatible pairs the structure of $G_{\gg,A}(\FF)$ is rather non-trivial even 
for classical Lie algebras $\gg=o_n(\kk)$ and $g=sp_n(\kk)$ and $A=M_n(\kk)$.
To demonstrate this, we explicitly compute the  ``Cartan subgroup'' of  
$G_{\gg,A}(\FF)$ (Proposition \ref{pr:cartan for classical groups}) as follows. 
For the above classical compatible pair $(\gg,A)$ an invertible diagonal matrix
$D=diag(f_1,...,f_n)\in M_n(\FF)$ belongs to the  $G_{\gg,A}(\FF)$ if and only if
$$f_if_{n-i+1}-f_1f_n\in I_1=\FF[\FF,\FF]$$
for $i=1,\ldots, n$.

Our computation of the  ``Cartan subgroup'' for  $\gg=sl_2(\kk)$ and $A=M_n(\kk)$
is dramatically harder and constitutes the second main result of this paper 
(Theorem \ref{th:sl2 in sln}). More precisely,
let $\gg=sl_2(\kk)\subset End_\kk(V)=M_n(\kk)=A$, where $V=S^{n-1}(\kk^2)$ is 
the simple $n$-dimensional $sl_2(\kk)$-module. Then an invertible diagonal matrix
$diag (f_1,f_2,\dots ,f_n)$ belongs to $G_{\gg,A}(\FF)$ if and only if
\begin{equation}
\label{eq:mainintro2}
\sum_{i=0}^m(-1)^i\binom {m}{i}f_{i+1}f_{i+2}^{-1}\in I_m
\end{equation}
for $m=1,2,\dots, n-2$.

We can also apply  our functorial generalization of $GL_n(\FF)$ to $K$-theory 
(however, we postpone all computations for concrete rings $\FF$ until a separate paper). We propose to generalize the
fundamental inclusion which plays a pivotal role in the algebraic
$K$-theory
\begin{equation}
\label{eq:pivotal}
E_n(\FF)\subset GL_n(\FF)
\end{equation}
where $E_n(\FF)$ is the subgroup generated by matrices
$1+tE_{ij},t\in \FF,  i\neq j$. Here the $E_{ij}$ are elementary
matrices. It is well-known and widely used that $E_n(\FF)$ is
normal in $GL_n(\FF)$ for $n\ge 3$ (and for certain
algebras $\FF $ when $n=2$).

It turns out that both groups $E_n(\FF)$ and $GL_n(\FF)$ are completely determined by the compatible pair $(\gg,A)$, where $\gg:=sl_n(\kk)\subset M_n(\kk)=A$.
Then $sl_n(\FF)=(\gg,A)(\FF)$, $GL_n(\FF)=G_{\gg,A}(\FF)$ and the group $E_n(\FF)$  is  generated by all elements $g\in G_{\gg,A}(\FF)$ of the form
$g=1+t\otimes E$
where $E\in \gg$ is a nilpotent in $A$, $t\in \FF$.

Motivated by this observation, we propose to generalize the
inclusion \eqref{eq:pivotal} to all of our $\NN$-current  groups  as follows.
%\comment{POVTORENIE?}Contrary to the ``straightforward" approach to classical groups over
%(noncommutative) algebras, we start this paper by
%introducing ``noncommutative'' Lie algebras as
%{\it noncommutative current  Lie algebras}, or $\NN$-current Lie algebras. For any pair $(\gg, A)$, where
%$\gg$ is a Lie  subalgebras of the associative algebra $A$ and any algebra $\FF$ we define
%$\NN$-current Lie algebra
%$(\gg, A)(\FF)$ to be
%the Lie subalgebra  of $\FF \otimes A$ generated by
%$\FF \otimes\gg$. In
%other words, $(\gg, A)(\FF)$ can be viewed as an
%$\FF$-envelope of $\gg $ in $\FF \otimes A$.
%Important examples of $A$ include the endomorphism algebra of a
%$\gg$-module or  the universal enveloping algebra of $\gg$ but we do not restrict
%ourselves to special cases.

%We explicitly compute $(\gg ,A)(\FF)$ for a large class of
%algebras including semisimple and (generalized) Kac-Moody Lie algebras (Theorem
%\ref{th:perfect}). 
We define (Section \ref{From NN-Lie algebras to NN-groups and generalized K1-theories}) the group $E_{\gg,A}(\FF)$ to be the sub-group of $G_{\gg,A}(\FF)$ generated by all elements of the form $1+t\otimes s$ where $t\in \FF$ and $s\in \gg$ is a nilpotent in $A$. Clearly, $E_{\gg,A}(\FF)$  is a normal subgroup of $G_{\gg,A}(\FF)$ which plays the role of $E_n(\FF)$ in $GL_n(\FF)$. This defines a generalized $K_1$-functor
$$\FF\mapsto K_1^{\gg,A}(\FF)=G_{\gg,A}(\FF)/E_{\gg,A}(\FF) \ .$$ 
In Section \ref{From NN-Lie algebras to NN-groups and generalized K1-theories} we also construct other  generalized $K_1$-functors in which we replace $E_{\gg,A}(\FF)$
by  $E_{S,\gg, A}(\FF)$, the subgroup generated 
generated by $1+S$, where $S$ is a $G_{\gg,A}(\FF)$-invariant subset of nilpotents in $\FF\otimes A$. However,  
computation of the generalized $K_1$-functors is beyond the scope of the present paper and will be performed in a separate publication.

The paper is organized as follows.

$\bullet$ Section \ref{sect: commutator expansions} contains some preliminary results on ideals in associative algebras $\FF$ generated by
$k$-th commutator spaces of $\FF$. Several key results are based on the Jacobi-Leibniz type identity \eqref{eq:Leibniz Jacobi identity}.

$\bullet$ In Section \ref{sect:Noncommutative current Lie algebras} we introduce $\NN$-Lie algebras and their important subclass:
$\NN$-current Lie algebras
$(\gg, A)(\FF)$ over Lie algebras $\gg$. As our first examples, we describe algebras
$(\gg, A)(\FF)$ for all classical Lie algebras.

$\bullet$ Section \ref{sect:upper bounds} contains upper bounds for $\NN$-current Lie algebras.

$\bullet$ Section \ref{sect: perfect pairs} contains our main result for
Lie algebras $(\gg, A)(\FF)$: upper bounds for algebras $(\gg,A)(\FF)$ coincide with them for
a large class of compatible pairs $(\gg,A)$ including all such pairs for semisimple Lie algebras $\gg$.

$\bullet$ In Section \ref{sect: n-groups} we introduce affine $\NN$-groups and $\NN$-current groups, there relation with $\NN$-Lie algebras and $\NN$-current Lie algebras,
important classes of their normal subgroups similar to subgroups $E_n(\FF)$ and the corresponding $K_1$-functors.
We also consider useful examples of $\NN$-subgroups and their ``Cartan subgroups" attached to the standard
representations of classical Lie algebras $\gg$ and to various
representations of $\gg=sl_2(\kk)$ Our description of these subgroups is based on a new class of algebraic identities for noncommutative difference derivatives (Lemma \ref{le:partial decomposition}) which are of interest by
themselves.

The present paper continues our study of algebraic groups over
noncommutative rings and their representations started in \cite{BR}. Part
of our results was published in \cite{BR1}.
In the next paper we will focus on  ``reductive groups over noncommutative rings",
their geometric structure and representations.

Throughout the paper, ${\bf Alg}$ will denote the category which objects are associative algebras (not necessarily with $1$)
over a field $\kk$ of characteristic zero and morphisms are algebra homomorphisms; and $\NN$ will stand for a sub-category of ${\bf Alg}$.
Also ${\bf Alg}_1$ will denote that sub-category of ${\bf Alg}$ which objects are unital $\kk$-algebras
over $\kk$  and arrows are homomorphisms of unital algebras.

\vspace{0.2cm}

\noindent {\bf Acknowledgements}. The authors would like to thank M. Kapranov
for very useful discussions and encouragements during
the preparation of the manuscript and C. Reutenauer for explaining
an important Jacobi type identity for commutators. The authors are grateful to
Max-Planck-Institut f\"ur Mathematik for its hospitality and
generous support during the essential stage of the work.

%\section{Noncommutative current Lie algebras}

\section{Commutator expansions and identities}
\label{sect: commutator expansions}

Given an object $\FF \in {\bf Alg}$,
for $k\ge 0$ define the $k$-th {\it commutator space $\FF^{(k)}$}
of $\FF$ recursively as $\FF^{(0)}=\FF$,
$\FF^{(1)}=\FF'=[\FF,\FF]$, $\FF^{(2)}=\FF''=[\FF,\FF']$, $\dots
$, $\FF^{(k)}=[\FF,\FF^{(k-1)}]$, $\dots $, where for any subsets
$S_1,S_2$ of $\FF$ the notation $[S_1,S_2]$ stands for the linear
span of all commutators $[a,b]=ab-ba$, $a\in S_1$, $b\in S_2$.
In a similar fashion, for each subset $S\subset \FF$ define the subspaces $S^{(k)}\subset \FF^{(k)}$ by
$S^{(0)}=span(S)$, $S^{(k)}=[S,S^{(k-1)}]$ for $k>0$; and the subspace $S^{(\bullet)}$  in
$\FF$ by
\begin{equation}
\label{eq:abstract Lie algebra}
S^{(\bullet)}=\sum_{k\ge 0} S^{(k)}
\end{equation}

The following result is obvious.

\begin{lemma} For any $S\subset \FF$ the subspace $S^{(\bullet)}$
is the Lie subalgebra of $\FF$ generated by $S$.
\end{lemma}

Following \cite{K} and \cite{KR}, define the subspaces $I_k^\ell(\FF)$ by:
$$I_k^\ell(\FF)=\sum_{\lambda} \FF^{(\lambda_1)} \FF^{(\lambda_2)} \cdots \FF^{(\lambda_\ell)} \ ,$$
where the summation goes over all
$\lambda=(\lambda_1,\lambda_2,\ldots,\lambda_\ell)\in (\ZZ_{ \ge 0})^\ell$ such that $\sum\limits_{i=1}^\ell\lambda_i=k$. Denote also
\begin{equation}
\label{eq:Iklel} I_k^{\le \ell}(\FF):=\sum_{1\le \ell'\le \ell}
I_k^{\ell'}(\FF), I_k (\FF):=I_k^{< \infty} =\sum\limits_{\ell \ge
1} I_k^\ell(\FF) \ .
\end{equation}
Clearly, $\FF I_k^\ell(\FF), I_k^\ell(\FF) \FF \subset
I_k^{\ell+1}(\FF)$. Therefore, $I_k(\FF)$ is a two-sided ideal in
$\FF$. Taking into account that
$\FF^{(k)}\FF+\FF^{(k)}=\FF\FF^{(k)}+\FF^{(k)}$ for all $k$, it
is easy to see that  $I_k(\FF)=I_k^k(\FF)+\FF I_k^k(\FF)$ for all
$k$.

\begin{lemma}
\label{le:embedded Ikl}
For each $k,\ell\ge 1$ one has:

\noindent (a) $I_k^\ell(\FF)\subset I_{k-1}^\ell(\FF)$, $I_k^{\le \ell}(\FF)\subset I_{k-1}^{\le \ell}(\FF)$.

\noindent (b) $[\FF,I_{k-1}^\ell(\FF)]\subset I_k^\ell(\FF)$, $[\FF,I_{k-1}^{\le \ell}(\FF)]\subset I_k^{\le \ell}(\FF)$.

\noindent (c) $I_k^{\le \ell+1}(\FF)=\FF[\FF,I_{k-1}^{\le \ell}(\FF)]+[\FF,I_{k-1}^{\le \ell+1}(\FF)]$.
\end{lemma}

\begin{proof} To prove (a) and (b), we need the following obvious recursion for $I_k^{\ell}(\FF)$, $\ell >1$:
\begin{equation}
\label{eq:recursion Ikl}
I_k^{\ell}(\FF)=\sum_{i\ge 0} \FF^{(i)} I_{k-i}^{\ell-1}(\FF)
\end{equation}
(with the natural convention that $I_{k'}^{\ell'}(\FF)=0$ if
$k'<0$). Then we prove (a) by induction in $\ell$. If $\ell=1$,
the assertion becomes $\FF^{(k)}\subset \FF^{(k-1)}$. Iterating
this inclusion and using the inductive hypothesis, we obtain for
$\ell >1$
$$I_k^{\ell}(\FF)=\sum_{i\ge 0} \FF^{(i)} I_{k-i}^{\ell-1}(\FF)= \FF I_k^{\ell-1}(\FF)+
\sum_{i> 0} \FF^{(i)} I_{k-i}^{\ell-1}(\FF)\subset $$
$$\subset \FF I_{k-1}^{\ell-1}(\FF)+\sum_{i>0} \FF^{(i-1)} I_{k-i}^{\ell-1}(\FF)= \sum_{i\ge 0} \FF^{(i)} I_{k-i-1}^{\ell-1}(\FF)=I_{k-1}^{\ell}(\FF) \ .$$
This proves (a).

Prove (b) also by induction in $\ell$. If $\ell=1$, the assertion becomes $[\FF,\FF^{(k-1)}]\subset \FF^{(k)}$, which is obvious. Using the inductive hypothesis, we obtain
$$[\FF,I_{k-1}^{\ell}(\FF)]=\sum_{i\ge 1} [\FF,\FF^{(i-1)} I_{k-i}^{\ell-1}(\FF)]\subset$$
$$\subset \sum_{i\ge 1} [\FF,\FF^{(i-1)}] I_{k-i}^{\ell-1}(\FF) + \FF^{(i-1)} [\FF,I_{k-i}^{\ell-1}(\FF)]
\subset \sum_{i\ge 0} \FF^{(i)} I_{k-i}^{\ell-1}(\FF)=I_k^{\ell}(\FF) \ .$$
This proves (b).

Prove (c). Obviously, $I_k^{\le \ell+1}(\FF)\supset \FF[\FF,I_{k-1}^{\le \ell}(\FF)]+[\FF,I_{k-1}^{\le \ell+1}(\FF)]$ by (b). Therefore, it suffices to prove the opposite inclusion
$$I_k^{\le \ell+1}(\FF)\subset \FF[\FF,I_{k-1}^{\le \ell}(\FF)]+[\FF,I_{k-1}^{\le \ell+1}(\FF)] \ .$$
We will use the following obvious consequence of \eqref{eq:recursion Ikl}:
$$I_k^{\le \ell+1}(\FF)=\sum_{i\ge 0} \FF^{(i)} I_{k-i}^{\le \ell}(\FF) \ .$$
Therefore, it suffices to prove that
\begin{equation}
\label{eq:hypothesis monomial}
\FF^{(i)} I_{k-i}^{\le \ell}(\FF)\subset \FF[\FF,I_{k-1}^{\le \ell}(\FF)]+[\FF,I_{k-1}^{\le \ell+1}(\FF)]
\end{equation}
for all $i\ge 0$, $\ell\ge 1$, $k\ge 1$. We prove  \eqref{eq:hypothesis monomial} by induction in all pairs $(\ell,i)$ 
ordered lexicographically. Indeed, suppose that the assertion is proved for all $(\ell',i')<(\ell,i)$. 
The base of induction is when $\ell=1$, $i=0$. Indeed, $I^{\le 1}_k(\FF)= \FF^{(k)}$ for all $k$ and 
\eqref{eq:hypothesis monomial} becomes
$\FF \FF^{(k)}\subset \FF[\FF,\FF^{(k-1)}]+[\FF,I_{k-1}^{\le 2}(\FF)]$, which is obviously true since $[\FF,\FF^{(k-1)}]=\FF^{(k)}$.

If $\ell\ge 1$, $i>0$, we obtain, using the Leibniz rule, the following inclusion:
$$\FF^{(i)} I_{k-i}^{\le \ell}(\FF)=[\FF,\FF^{(i-1)}] I_{k-i}^{\le \ell}(\FF)\subset [\FF,\FF^{(i-1)} I_{k-i}^{\le \ell}(\FF)]+ \FF^{(i-1)} [\FF,I_{k-i}^{\le \ell}(\FF)] \ .$$
Therefore, by b)$$\FF^{(i)} I_{k-i}^{\le \ell}(\FF)\subset
[\FF,\FF^{(i-1)} I_{k-i}^{\le \ell}(\FF)]+ \FF^{(i-1)}
I_{k+1-i}^{\le \ell}(\FF).$$

Finally, using the inductive hypothesis for $(\ell,i-1)$ and taking into account that $\FF^{(i-1)} I_{k-i}^{\le \ell}(\FF)\subset I_{k-1}^{\le \ell+1}(\FF)$, and, therefore, $$[\FF,\FF^{(i-1)} I_{k-i}^{\le \ell}(\FF)]\subset [\FF,I_{k-1}^{\le \ell+1}(\FF)]\ ,$$
we obtain the inclusion \eqref{eq:hypothesis monomial}.

If $\ell\ge 2$, $i=0$,  then using the inductive hypothesis for all pairs $(\ell-1,i')$, $i'\ge 0$, we obtain:
 $$I_k^{\le \ell}(\FF)=\FF[\FF,I_{k-1}^{\le \ell-1}(\FF)]+[\FF,I_{k-1}^{\le \ell}(\FF)] \ .$$
Multiplying by $\FF$ on the left
%and using the distributivity of multiplication of subspaces in $\FF\cdot A$,
we obtain:
$$\FF I_k^{\le \ell}(\FF)=\FF^2[\FF,I_{k-1}^{\le \ell-1}(\FF)]+\FF[\FF,I_{k-1}^{\le \ell}(\FF)] =\FF[\FF,I_{k-1}^{\le \ell}(\FF)] $$
because $\FF^2\subset \FF$ and $I_{k-1}^{\le \ell-1}(\FF)\subset I_{k-1}^{\le \ell}(\FF)$. This immediately implies \eqref{eq:hypothesis monomial}.

Part (c) is proved.
The lemma is proved.
\end{proof}

\begin{lemma}
\label{le:ideals I_k(F) simple}
For any $k',k\ge 0$, and any $\ell,\ell'\ge 1$ one has:

\noindent (a) $I_k^\ell(\FF) I_{k'}^{\ell'}(\FF)\subset I_{k+k'}^{\ell+\ell'}(\FF)$, $I_k^{\le \ell}(\FF) I_{k'}^{\le \ell'}(\FF)\subset I_{k+k'}^{\le \ell+\ell'}(\FF)$.

\noindent (b) $[I_k^\ell(\FF), I_{k'}^{\ell'}(\FF)]\subset [\FF,I_{k+k'}^{\ell+\ell'-1}(\FF)]$, $[I_k^{\le \ell}(\FF), I_{k'}^{\le \ell'}(\FF)]\subset [\FF,I_{k+k'}^{\le \ell+\ell'-1}(\FF)]$.

\end{lemma}

\begin{proof} Part (a) follows from the obvious fact that
$$(\FF^{(\lambda_1)} \FF^{(\lambda_2)} \cdots \FF^{(\lambda_{\ell_1})} )
(\FF^{(\mu_1)} \FF^{(\mu_2)} \cdots \FF^{(\mu_{\ell_2})})\subset I_k^{\ell_1+\ell_2}(\FF)\ ,$$
where $k=\lambda_1+\lambda_2+\cdots+\lambda_{\ell_1}+\mu_1+\mu_2+\cdots+\mu_{\ell_2}$.

Prove (b).
First, we prove the first inclusion for $\ell=1$. We proceed by induction on $k$. 
The base of induction, $k=0$, is obvious because $I^1_0(\FF)=\FF$. Assume that the assertion is proved for all $k_1<k$, i.e., we have:
$$[\FF^{(k_1)}, I_{k'}^{\ell'}(\FF)]\subset [\FF,I_{k_1+k'}^{\ell'}(\FF)] \ .$$
Then, using the fact that $\FF^{(k)}=[\FF,\FF^{(k-1)}]$ and the Jacobi identity, we obtain:
$$[\FF^{(k)}, I_{k'}^{\ell'}(\FF)]=[[\FF,\FF^{(k-1)}], I_{k'}^{\ell'}(\FF)]\subset$$
$$\subset [\FF,[\FF^{(k-1)}, I_{k'}^{\ell'}(\FF)]]+[\FF^{(k-1)}, [\FF,I_{k'}^{\ell'}(\FF)]]\subset $$
$$\subset [\FF, [\FF,I_{k'+k-1}^{\ell'}(\FF)]]+[\FF^{(k-1)}, I_{k'+1}^{\ell'}(\FF)]\subset [\FF,I_{k'+k}^{\ell'}(\FF)] $$
by the inductive hypothesis and Lemma \ref{le:embedded Ikl}(b). This proves the first inclusion of (b) for $\ell=1$.

Furthermore, we will proceed by induction on $\ell$. Now $\ell>1$, assume that the assertion is proved for all $\ell_1<\ell$, i.e., we 
have the inductive hypothesis in the form:
$$[I_k^{\ell_1}(\FF), I_{k'}^{\ell'}(\FF)]\subset [\FF,I_{k+k'}^{\ell_1+\ell'-1}(\FF)]$$
for all $k, k'\ge 0$.

We need the following useful Jacobi-Leibniz type identity in
$\FF$:
\begin{equation}
\label{eq:Leibniz Jacobi identity}
[ab,c]+[bc,a]+[ca,b]=0
\end{equation}
for all $a,b,c\in \FF$ (The identity was communicated to the authors by C. Reutenauer and was
used in a different context in the recent paper \cite{FeSh}).

Using  \eqref{eq:recursion Ikl} and \eqref{eq:Leibniz Jacobi identity} with all $a\in \FF^{(i)}$, $b\in I_{k-i}^{\ell-1}(\FF)$, $c\in I_{k'}^{\ell'}(\FF)$, we obtain for all $i\ge 0$:
$$[\FF^{(i)} I_{k-i}^{\ell-1}(\FF), I_{k'}^{\ell'}(\FF)]\subset [\FF^{(i)}, I_{k-i}^{\ell-1}(\FF) I_{k'}^{\ell'}(\FF)]+
[I_{k-i}^{\ell-1}(\FF), I_{k'}^{\ell'}(\FF)\FF^{(i)}]\subset $$
$$[\FF^{(i)}, I_{k+k'-i}^{\ell+\ell'-1}(\FF)]+[I_{k-i}^{\ell-1}(\FF),I_{k'+i}^{\ell'+1}(\FF)]\subset [\FF,I_{k+k'}^{\ell+\ell'-1}(\FF)]$$
by the the already proved (a) and inductive hypothesis. This finishes the proof of the first inclusion of (b). The second inclusion of (b) also follows.
\end{proof}

Generalizing \eqref{eq:recursion Ikl}, for any subset $S$ of
$\FF$ denote by $I_k^\ell(\FF,S)$ the image of
$Span\ S^{\otimes(k+\ell)}$ under the canonical map $\FF^{\otimes
(k+\ell)}\twoheadrightarrow I_k^\ell(\FF)$, i.e.,
\begin{equation}
\label{eq:Ikl(S)} I_k^\ell(\FF,S)=\sum_{\lambda} S^{(\lambda_1)}
S^{(\lambda_2)} \cdots S^{(\lambda_\ell)} \ ,
\end{equation}
$\lambda=(\lambda_1,\lambda_2,\ldots,\lambda_\ell)\in (\ZZ_{ \ge
0})^\ell$ such that $\sum\limits_{i=1}^\ell\lambda_i=k$.

In particular, $I_k^1(\FF,S)=S^{(k)}$ and $I_0^\ell=S^\ell$.

The following result is obvious.

\begin{lemma}
\label{le:ideals I_k}  Let $\FF$ be an object of $\NN$ and $S\subset \FF$. Then:

\noindent (a)
For any $k\ge 0$, $\ell\ge 2$ one has
$$I_k^\ell(\FF,S)=\sum_{i=0}^k S^{(i)} I_{k-i}^{\ell-1}(\FF,S) \ .$$
\noindent (b) For any $k',k\ge 0$, and any $\ell,\ell'\ge 1$ one has:
$$I_k^\ell(\FF,S) I_{k'}^{\ell'}(\FF,S)\subset I_{k+k'}^{\ell+\ell'}(\FF, S),~[I_k^\ell(\FF,S), I_{k'}^{\ell'}(\FF,S)]\subset I_{k+k'+1}^{\ell+\ell'-1}(\FF, S)\ .$$
%\noindent (iii) If $k<\ell$ and $S'\subset S$, then $I_\ell\subset I_k$.
In particular,
\begin{equation}
\label{eq:ideals I_k simple}
S^{(i)}I_k^\ell(\FF,S)\subset I_{k+i}^{\ell+1}(\FF,S),~[S^{(i)},I_k^\ell(\FF,S)]\subset I_{k+i+1}^\ell(\FF,S) \ .
\end{equation}
\end{lemma}

\section{$\NN$-Lie algebras and $\NN$-current Lie algebras}
\label{sect:Noncommutative current Lie algebras}

Given objects $\FF$  and ${\mathcal A}$ of ${\bf Alg}$, we refer to a morphism $\iota:\FF\to {\mathcal A}$ in  
${\bf Alg}$ as an  {\it $\FF$-algebra structure on ${\mathcal A}$} (we will also refer to ${\mathcal A}$ as an {\it  $\FF$-algebra}).

Note that each $\FF$-algebra structure on ${\mathcal A}$ turns ${\mathcal A}$ into an algebra  in the category of $\FF$-bimodules
(i.e.,   ${\mathcal A}$ admits two $\FF$-actions $\FF\otimes {\mathcal A}\to {\mathcal A}$, ${\mathcal A}\otimes \FF\to {\mathcal A}$ 
via  $f\otimes a\mapsto \iota(f)\cdot a$ and $a\otimes f\mapsto a\cdot \iota(f)$ respectively).

We fix an arbitrary sub-category $\NN$ of ${\bf Alg}$ throughout
the section. In most cases we take $\NN={\bf Alg}$.

\begin{definition}
\label{def: FF-Lie algebra} An $\NN$-Lie algebra is a triple $(\FF,{\mathcal L},{\mathcal A})$, where 
$\FF$ is an object of $\NN$, ${\mathcal A}$ is an $\FF$-algebra, 
and  ${\mathcal L}$ is an $\FF$-{\it Lie subalgebra of ${\mathcal A}$}, i.e.,  
if ${\mathcal L}$ is a Lie subalgebra (under the commutator bracket) of ${\mathcal A}$  
invariant under the adjoint action of $\FF$ on ${\mathcal A}$ given by 
$(f,a)\mapsto \iota(f)\cdot a-a\cdot \iota(f)$ for all $f\in \FF$, $a\in {\mathcal A}$.
\end{definition}

A morphism $(\FF_1,{\mathcal L}_1,{\mathcal A}_1)\to (\FF_2,{\mathcal L}_2,{\mathcal A}_2)$ of 
$\NN$-Lie algebras is a pair $(\varphi,\psi)$, where $\varphi:\FF_1\to \FF_2$ 
is a morphism in $\NN$ and  $\psi:{\mathcal A}_1\to {\mathcal A}_2$ 
is a morphism in ${\bf Alg}$ such that  $\psi({\mathcal L}_1)\subset {\mathcal L}_2$ and
$\psi\circ \iota_1=\iota_2\circ \varphi$.

Denote by ${\bf LieAlg}_\NN$ the category of $\NN$-Lie algebras.

For an $\NN$-Lie algebra $(\FF,{\mathcal L},{\mathcal A})$, let
$L_i(\FF,{\mathcal L},{\mathcal A}):=(\FF,{\mathcal L}^{(i)},{\mathcal A})$, $0\le i\le j\le 3$, 
where ${\mathcal L}^{(i)}$, $i=1,2,3$ are given by:

\noindent $\bullet$ ${\mathcal L}^{(1)}$ is the normalizer Lie algebra of ${\mathcal L}$ in ${\mathcal A}$.

\noindent $\bullet$ ${\mathcal L}^{(2)}$ is the Lie subalgebra of ${\mathcal A}$ 
generated by $\iota(\FF)\subset {\mathcal A}$ and by the semigroup
${\mathcal S}=\{s\in {\mathcal A}:~s\cdot {\mathcal L}={\mathcal L}\cdot s\}$.

\noindent $\bullet$ ${\mathcal L}^{(3)}$ is the Lie subalgebra of ${\mathcal A}$ generated by 
$\GG(\iota(\FF))\subset {\mathcal A}$, where ${\mathcal G}$ is the stabilizer of 
${\mathcal L}$ in the group $Aut_\kk({\mathcal A})$, i.e.,
\begin{equation}
\label{eq:NN-Lie algebras}{\mathcal G}=\{g\in Aut_\kk({\mathcal A}): g({\mathcal L})={\mathcal L}\} \ .
\end{equation}

The following result is obvious.

\begin{lemma}
\label{le:Lie functors} For each $\NN$-Lie algebra $(\FF,{\mathcal L},{\mathcal A})$ and $i=1,2,3$ the triple
$L_i(\FF,{\mathcal L},{\mathcal A})$ is also an $\NN$-Lie algebra.

\end{lemma}

Therefore, we can construct a number of new $\NN$-Lie algebras by combining the operations $L_i$ for a given $\NN$-Lie algebra. 

\begin{remark} In general, none of $L_i$ defines a functor ${\bf LieAlg}_\NN\to {\bf LieAlg}_\NN$. 
However, for each $i=1,2,3$ one can find an appropriate subcategory ${\mathcal C}$ 
of ${\bf LieAlg}_\NN$ such the restriction of $L_i$ to ${\mathcal C}$ is a functor ${\mathcal C}\to {\bf LieAlg}_\NN$.
\end{remark}

\begin{remark} The operation $L_3$ is interesting only when $\FF$ is noncommutative because 
for any object  $(\FF,{\mathcal L},{\mathcal A})$ of ${\bf LieAlg}_\NN$ such that $\FF$ is 
commutative and all automorphisms of ${\mathcal A}$ are inner, one obtains ${\mathcal L}^{(3)}=\iota(\FF)$ 
and therefore, $L_3(\FF,{\mathcal L},{\mathcal A})=(\FF,\iota(\FF),{\mathcal A})$.

\end{remark}

Denote by  $\pi$  the natural (forgetful) projection functor ${\bf
LieAlg}_\NN\to \NN$ such that $\pi(\FF,{\mathcal L},{\mathcal A})=\FF$ and $\pi(\varphi,\psi)=\varphi$.
\begin{definition} A {\it noncommutative current Lie algebra} ({\it $\NN$-current Lie algebra}) is a functor ${\mathfrak s}:\NN\to
{\bf
LieAlg}_\NN$ such that $\pi\circ {\mathfrak s}=Id_{\NN}$ (i.e., ${\mathfrak s}$ is a section of $\pi$).
\end{definition}

Note that if $\NN=(\FF,Id_\FF)$ has only one object $\FF$ and only the identity arrow $Id_\FF$, 
then the $\NN$-current Lie algebra is simply any object of ${\bf LieAlg}_{\bf Alg}$ of the
form $(\FF,{\mathcal L},{\mathcal A})$. In this case, we will sometimes refer to the Lie algebra ${\mathcal L}$ as a {\it 
$\FF$-current Lie algebra}.

In principle, we can construct a number of $\FF$-current Lie algebras by twisting a given one with operations $L_i$ 
from Lemma \ref{le:Lie functors}. However, the study of such ``derived'' $\FF$-current Lie algebras 
is beyond the scope of the present paper.

In what follows we will suppress the tensor sign in expressions like $\FF\otimes A$ and 
write $\FF\cdot A$ instead. Note that for any object $A$ of ${\bf Alg}_1$ and 
any object $\FF$ of ${\bf Alg}$ the product $\FF\otimes A$ is naturally an $\FF$-algebra 
via the embedding $\FF\hookrightarrow \FF\cdot A$ ($f\mapsto f\cdot 1$).

The following is a first obvious example of $\NN$-current Lie algebras.

\begin{lemma} For any object algebra $A$ of ${\bf Alg}_1$ and any object $\FF$ of $\NN$ define the object ${\mathfrak
s}_A(\FF)=(\FF,\FF\cdot A,\FF\cdot A)$ of ${\bf LieAlg}_\NN$. Then the  association $\FF\mapsto {\mathfrak s}_A(\FF)$ defines a
noncommutative current Lie algebra ${\mathfrak s}_A:\NN\to {\bf LieAlg}_\NN$.
\end{lemma}

The main object of our study will be a refinement of the above example.
Given an object $A$ of ${\bf Alg}_1$, and a subspace $\gg\subset
A$ such that $[\gg,\gg]\subset \gg$ (i.e., $\gg$
is a Lie subalgebra of $A$), we say that $(\gg,A)$ is a {\it compatible pair}.
For any compatible pair $(\gg,A)$ and  an  object $\FF$ of $\NN$, denote by  $(\gg,A)(\FF)$ the Lie subalgebra of the
$\FF\cdot A=\FF\otimes A$ (under the commutator bracket) generated by $\FF\cdot \gg$, that is,
$(\gg,A)(\FF)=(\FF\cdot \gg)^{(\bullet)}$ in notation \eqref{eq:abstract Lie algebra}.

\begin{proposition} For any compatible pair $(\gg, A)$ the association
$$\FF\mapsto (\FF,(\gg,A)(\FF),\FF\cdot A)$$ defines the $\NN$-current Lie algebra
$$(\gg,A):\NN\to {\bf LieAlg}_\NN \ .$$

\end{proposition}

\begin{proof} It suffices to show that any arrow $\varphi$ in $\NN$, i.e., any algebra 
homomorphism $\varphi: \FF_1\to \FF_2$ defines a homomorphism of Lie algebras 
$(\gg,A)(\FF_1)\to (\gg,A)(\FF_2)$. We need the following obvious fact.

\begin{lemma}
\label{le:NA-functoriality}
Let ${\mathcal A}_1$, ${\mathcal A}_2$ be objects of ${\bf Alg}$ and let $\varphi:A_1\to A_2$ be 
a morphism in ${\bf Alg}$. Let $S_1\subset {\mathcal A}_1$ and 
$S_2\subset {\mathcal A}_2$ be two subsets such that $\varphi(S_1)\subset S_2$. 
Then the restriction of $\varphi$ to the Lie algebra $S_1^{(\bullet)}$ 
(in notation \eqref{eq:abstract Lie algebra}) is a homomorphism of Lie algebras $S_1^{(\bullet)}\to S_2^{(\bullet)}$.
\end{lemma}

Applying Lemma \ref{le:NA-functoriality} with ${\mathcal A}_i=\FF_i\cdot A$, $S_i=\FF_i\cdot \gg$, $i=1,2$, $\varphi=f\otimes 
id_A:\FF_1\cdot A\to \FF_2\cdot A$, 
%the trivial extension of $\FF$, 
we obtain a Lie algebra homomorphism $(\gg,A)(\FF_1)=(\FF_1\cdot \gg)^{(\bullet)}\to  (\FF_2\cdot 
\gg)^{(\bullet)}=(\gg,A)(\FF_2)$.

It remains to show that the action of $\FF$ on ${\mathcal
L}=(\gg,A)(\FF)=(\FF\cdot \gg)^{(\bullet)}$ is stable under the
commutator bracket with $\FF$. Indeed, $S=\FF\cdot \gg$ is
invariant under the adjoint action of $\FF$ on $\FF\cdot A$. By
induction and the Jacobi identity ${\mathcal L}=(\FF\cdot
\gg)^{(\bullet)}$ is also invariant under this action of $\FF$.
The proposition is proved.
\end{proof}

If $\FF$ is commutative, then $(\gg,A)(\FF)=\FF\cdot \gg$ is the
$\FF$-current algebra. Therefore, if $\FF$ is an arbitrary object
of $\NN$, the Lie algebra $(\gg,A)(\FF)$ deserves a name of the
{\it $\NN$-current Lie algebra associated with the compatible pair $(\gg,A)$}.

%We will refer to this functor $(\gg,A)$ as the noncommutative current Lie algebra {\it associated with the compatible pair $(\gg,A)$}.

If $A=U(\gg)$, the universal enveloping algebra of $\gg$, then we abbreviate $\gg(\FF):=(\gg,U(\gg))(\FF)$. 
Another natural choice of $A$ is algebra $End(V)$, where $V$ is a faithful $\gg$-module. 
In this case, we will sometimes abbreviate $(\gg,V)(\FF):=(\gg,End(V))(\FF)$.

The following result provides an estimation of $(\gg,A)(\FF)$ from below. Set
\begin{equation}
\label{eq:langle gg rangle}
\langle \gg\rangle=\sum_{k\ge 1} \gg^k \ ,
\end{equation}
i.e.,   $\langle \gg\rangle$ is the associative subalgebra of $A$
%, maybe without unit, 
generated by $\gg$.

\begin{proposition}
\label{pr:lower bound small}
Let $(\gg,A)$ be a compatible pair and $\FF$ be an object of $\NN$. Then

\noindent (a) $\FF^{(k)}\cdot \gg^{k+1}\subset (\gg,A)(\FF)$ and $\FF \FF^{(k)}
\cdot [\gg,\gg^{k+1}]\subset (\gg,A)(\FF)$ for all $k\ge 0$.

\noindent (b) If $\gg$ is abelian, i.e., $\gg'=[\gg,\gg]=0$, then
\begin{equation}
\label{eq:abelian currents}(\gg,A)(\FF)=\sum_{k\ge 0} \FF^{(k)}\cdot \gg^{k+1} \ .
\end{equation}

\noindent (c) If $[\FF,\FF]=\FF$ (i.e., $\FF$ is perfect as a Lie algebra), then
$(\gg,A)(\FF)=\FF\cdot \langle \gg\rangle$.
\end{proposition}

\begin{proof} Prove (a). We need  the following technical result.

\begin{lemma}
\label{le:powers of gg}
Let $(\gg, A)$ be a compatible pair.
For all $m\ge 2$ denote by $\widetilde {\gg^m}$ the $\kk$-linear span of all powers $g^m$, $g \in \gg$. Then for any $m\ge 2$ one has
\begin{equation}
\label{eq:tilde 1} \widetilde {\gg^m}+(\gg^{m-1}\cap \gg^m)=\gg^m
\end{equation}
\end{lemma}

\begin{proof} Since $\gg^{i-1}\gg'\gg^{m-i-1}\subset \gg^{m-1}$ for all $i\le m-1$, we obtain the following congruence for any
$c=(c_1,\ldots,c_m)$, $c_i\in (\kk-\{0\})$ and
$x=(x_1,\ldots,x_m)$, $x_i \in \gg$, $i=1,2,\dots ,m$:
$$(c_1x_1+\cdots c_m x_m)^m\equiv \sum_{\lambda} \binom{m}{\lambda} c^\lambda x^\lambda \mod (\gg^{m-1}\cap \gg^m)\ ,$$
where the summation is over all partitions
$\lambda=(\lambda_1,\ldots,\lambda_m)$ of $m$ and we abbreviated
$c^\lambda=c_1^{\lambda_1}\cdots c_m^{\lambda_m}$ and
$x^\lambda=x_1^{\lambda_1}\cdots x_m^{\lambda_m}$. Varying
$c=(c_1,\cdots,c_m)$, the above congruence implies that that each
monomial $x^\lambda$ belongs to $\widetilde
{\gg^m}+(\gg^{m-1}\cap \gg^m)$. In particular, taking
$\lambda=(1,1,\ldots,1)$, we obtain $\gg^m\subseteq \widetilde
{\gg^m}+(\gg^{m-1}\cap \gg^m)$. Taking into account that
$\widetilde {\gg^m}\subseteq \gg^m$, we obtain \eqref{eq:tilde 1}.
The lemma is proved.
\end{proof}

We also need the following useful identity in $\FF\cdot A$:
\begin{equation}
\label{eq:the commutator noncom coeffs}
[sE,tF]=st\cdot [E,F]+[s,t]\cdot FE=ts\cdot [E,F]+[s,t]\cdot EF
\end{equation}
for any $s,t\in \FF$, $E,F\in A$.

We prove the first inclusion of (a)  by induction on $k$. If
$k=0$, one obviously has $\FF^{(0)}\gg^1=\FF\cdot \gg\subset
(\gg,A)(\FF)$. Assume now that $k>0$. Then for $g\in \gg$ we
obtain by using \eqref{eq:the commutator noncom coeffs}:
$$[\FF\cdot g,\FF^{(k-1)}\cdot g^k]=[\FF,\FF^{(k-1)}]\cdot g^{k+1}=\FF^{(k)}\cdot g^{k+1}$$
which implies that $\FF^{(k)}\cdot \widetilde {\gg^{k+1}}\subset (\gg,A)(\FF)$ 
(in the notation of Lemma \ref{le:powers of gg}). Using Lemma \ref{le:powers of gg}, we obtain
$$\FF^{(k)}\cdot \widetilde {\gg^{k+1}}\equiv \FF^{(k)}\cdot \gg^{k+1}\mod \FF^{(k)}\cdot (\gg^k\cap \gg^{k+1}) \ .$$
Taking into account that $\FF^{(k)}\cdot (\gg^k\cap
\gg^{k+1})\subset \FF^{(k-1)}\cdot \gg^k \subset (\gg,A)(\FF)$ by
the inductive hypothesis (here we used the inclusion
$\FF^{(k)}\subset \FF^{(k-1)}$), the above relation implies that
$\FF^{(k)}\cdot \gg^{k+1}$ also belongs to $(\gg,A)(\FF)$. This
proves the first inclusion of (a). To prove the second inclusion,
we compute using \eqref{eq:the commutator noncom coeffs}
$$[\FF\cdot \gg,\FF^{(k-1)}\cdot \gg^k]\equiv \FF \FF^{(k-1)}\cdot [\gg,\gg^k] \mod \FF^{(k)}\cdot \gg^{k+1} \ .$$
Therefore, using the already proved inclusion $\FF^{(k)}\cdot \gg^{k+1}\subset (\gg,A)(\FF)$, 
we see that $\FF \FF^{(k-1)}\cdot [\gg,\gg^k]$ also belongs to $(\gg,A)(\FF)$. This finishes the proof of  (a).

Prove (b). Clearly, (a) implies that $(\gg,A)(\FF)$ contains the
right hand side of \eqref{eq:abelian currents}. Therefore, it
suffices to prove that the latter space is closed under the
commutator. Indeed, since $\gg$ is abelian, one has
$$[\FF^{(k_1)}\cdot \gg^{k_1+1}, \FF^{(k_2)}\cdot \gg^{k_2+1}]=[\FF^{(k_1)},\FF^{(k_2)}]\cdot \gg^{k_1+k_2+2}\subset $$
$$\subset \FF^{(k_1+k_2+1)}\cdot \gg^{k_1+k_2+2}\subset (\gg,A)(\FF)$$
because $[\FF^{(k_1)},\FF^{(k_2)}]\subset \FF^{(k_1+k_2+1)}$.
This finishes the proof of (b).

Prove (c). Since $\FF'=\FF$, the already proved part (a) implies $\FF\cdot \gg^k\subset (\gg,A)(\FF)$ for all $k\ge 1$, 
therefore, $\FF\cdot \langle \gg \rangle \subseteq (\gg,A)(\FF)$. Since $\langle \gg \rangle$ is an associative subalgebra of $A$ 
containing $\gg$, we obtain an opposite inclusion $(\gg,A)(\FF)\subseteq \FF\cdot \langle \gg \rangle$. This finishes the proof of (c).

The proposition is proved.
\end{proof}

\begin{remark} Proposition \ref{pr:lower bound small}(c) shows that the case 
when $[\FF,\FF]=\FF$ is not of much interest. This happens, for example, when 
$\FF$ is a Weyl algebra or the quantum torus. In these cases a natural anti-involution on 
$\FF$ can be taken into account. We will discuss it in a separate paper.

\end{remark}
%Denote by $\langle \gg\rangle$ the subalgebra of $A$ generated by $\gg$.

%%%%%%%%%%%%%%%%%%%%%%%%%%%
%Let $\gg\subset gl(V)$ be a Lie subalgebra.
%Denote by $\pi:End(V)\to sl(V)$ the projection to the first summand of the canonical decomposition $gl(V)=sl(V)+V_0$, where $V_0=\kk\cdot 1$ is the (one-dimensional) center of $gl(V)$. For each $k\ge 1$ denote
%$$V_k:=\pi(\gg^k) \ ,$$
%where the power is taken in $gl(V)$ considered as the associative algebra $End(V)$. That is, $V_k=\{x\in %\gg^k:Tr(x)=0\}$.

\begin{definition}
\label{def:type} We say that a compatible pair $(\gg,A)$ is of
{\it finite type} if there exist $m>0$  such that
$\gg+\gg^2+\cdots +\gg^m=A$, and we call such minimal $m$ the
{\it type} of $(\gg,A)$. If such $m$ does not exists, we say that
$(\gg,A)$ is of infinite type.
\end{definition}
Note that $(\gg,A)$ is of type $1$ if and only if $\gg=A$, which, in its turn, implies that $(\gg,A)(\FF)=\FF\cdot A$ for all objects $\FF$ of $\NN$. Note also that if $\langle g\rangle =A$ and $A$ is finite-dimensional over $\kk$, then $(\gg,A)$ is always of finite type.

\begin{proposition}
\label{pr:type two} Assume that $(\gg,A)$ is of type $2$, i.e., $\gg\ne A$ and  $\gg+\gg^2=A$. Then
\begin{equation}
\label{eq:type two} (\gg,A)(\FF)=\FF\cdot \gg+\FF'\cdot A +\FF\FF' \cdot [A,A] \ ,
\end{equation}
where $\FF'=[\FF,\FF]$.
\end{proposition}

\begin{proof}
%Then we obtain the following congruence
%$$(\FF\cdot \gg)^{(1)}=[\FF\cdot \gg,\FF\cdot \gg]\equiv \FF'\cdot \gg^2 \mod \FF^2[\gg,\gg] \ .$$
%Since $\FF^2[\gg,\gg]\subset \FF\cdot \gg\subset (\gg,A)(\FF)$, we obtain $\FF'\cdot \gg^2\subset (\gg,A)(\FF)$. 
Since %$\FF'\subset \FF$, we obtain $\FF'\cdot (\gg+\gg^2)=\FF'\cdot A\subset  (\gg,A)(\FF)$. 
Furthermore,
%$$[\FF\cdot \gg,\FF'\cdot A]\equiv \FF\FF'\cdot [\gg,A] \mod [\FF,\FF'] A\gg \ .$$
%Since $[\FF,\FF'] A\gg\subset \FF'\cdot A\subset (\gg,A)(\FF)$, we obtain $\FF\FF'\cdot [\gg,A] \subset (\gg,A)(\FF)$.
%
Proposition \ref{pr:lower bound small}(a) guarantees that
$$\FF\cdot \gg+\FF'\cdot \gg^2 +\FF\FF' \cdot [\gg,\gg^2]\subset (\gg,A)(\FF) \ .$$
Clearly, $\FF\cdot \gg+\FF'\cdot \gg^2=\FF\cdot \gg+\FF'\cdot A$ (because $\FF'\subset \FF$). 
Let us now prove that $[\gg,A]=[A,A]$. Obviously, $[\gg,A]\subseteq [A,A]$. 
The opposite inclusion immediately follows from inclusion $[\gg^2,\gg^2]\subseteq [\gg,\gg^3]$, which, in its turn follows from 
\eqref{eq:Leibniz Jacobi identity}: taking any $a\in \gg$, $b\in \gg$, $c\in \gg^2$ in \eqref{eq:Leibniz Jacobi identity}, 
we obtain $[ab,c]\in [\gg,\gg^3]$.

Using the equation $[\gg,A]=[A,A]$ we obtain $\FF\FF'\cdot [A,A] \subset (\gg,A)(\FF)$. 
This proves that $(\gg,A)(\FF)$ contains the right hand side of \eqref{eq:type two}.

To finish the proof, it suffices to show that the latter subspace is closed under the commutator.
Indeed,  abbreviating $A'=[A,A]$, we obtain
$$[\FF\cdot \gg,\FF'\cdot A]\subset \FF\FF'\cdot [\gg,A] +[\FF,\FF']\cdot A\gg\subset \FF\FF'\cdot A'+\FF'\cdot A \ ,$$
$$[\FF\cdot \gg,\FF\FF'\cdot A']\subset \FF^2\FF'\cdot [\gg,A'] +[\FF,\FF\FF']\cdot A'\gg\subset \FF\FF'\cdot A'+\FF'\cdot A \ ,$$
$$[\FF'\cdot A,\FF'\cdot A]\subset (\FF')^2\cdot A' +[\FF',\FF']\cdot A^2\subset \FF\FF'\cdot A'+\FF'\cdot A \ ,$$
$$[\FF'\cdot A,\FF\FF' \cdot A']\subset \FF'\FF\FF'\cdot [A,A']+[\FF',\FF \FF']\cdot A'A  $$
because 
$\FF'\FF\FF'\cdot [A,A']\subset \FF\FF'\cdot A'\subset (\gg,A)(\FF)$ and 
$[\FF',\FF \FF']\cdot A'A\subset \FF'\cdot A\subset (\gg,A)(\FF)$. Finally,
$$[\FF'\FF\cdot A',\FF\FF' \cdot A']\subset (\FF'\FF)^2\cdot [A',A']+[\FF'\FF,\FF \FF']\cdot (A')^2 $$
because $(\FF'\FF)^2\cdot [A',A']\subset \FF'\FF\cdot A'$ and
$[\FF'\FF,\FF \FF']\cdot [A,A]^2\subset \FF'\cdot A$. The
proposition is proved.
\end{proof}

For any $\kk$-vector space $V$ and any object $\FF$ of $\NN$ we
abbreviate $sl(V,\FF):=(sl(V),End(V))(\FF)$ and
$gl(V,\FF):=(End(V),End(V))(\FF)=\FF\cdot End(V)$. The following
result shows that for $n=\dim V\ge 2$ the commutator Lie algebra $sl_n(\FF)=[gl_n(\FF), gl_n(\FF)]$ is, in fact,  $sl(V,\FF)$.

\begin{corollary} Let $V$ be a finite-dimensional $\kk$-vector space such that $\dim V>1$. 
Then $(\gg,A)=(sl(V),End(V))$ is of type $2$ and
$$sl(V,\FF)=\FF'\cdot 1+\FF\cdot sl(V)\ .$$
\end{corollary}

Hence  $sl(V,\FF)$ is the set of all $X\in  gl(V,\FF)$ such that $Tr(X)\in \FF'=[\FF,\FF]$  (where $Tr:gl(V,\FF)=\FF\cdot
End(V)\to \FF$ is the trivial extension of the ordinary trace $End(V)\to \kk$).

\begin{proof} Let us prove that the pair $(\gg,A)=(sl(V),End(V))$ is of type $2$, i.e., $sl(V)+sl(V)^2=End(V)$.
It suffices to show that $1\in sl(V)^2$. To prove it, choose a basis $e_1,\ldots,e_n$ in $V$ so that
$V\cong \kk^n$,  $sl(V)\cong sl_n(\kk)$ and $A=End(V)\cong M_n(\kk)$.
Indeed, for any indices $i\ne j$ both $E_{ij}$ and $E_{ji}$ belong to $sl(V)$,
and $E_{ij}E_{ji}=E_{ii}\in sl(V)^2$. Therefore, $1=\sum_{i=1}^n E_{ii}$ also belongs to $sl(V)^2$.
Applying Proposition \ref{pr:type two} and using the obvious fact that $[A,A]= sl(V)$, we obtain
$$sl(V,\FF)=\FF\cdot sl(V)+\FF'\cdot A+\FF\FF'[A, A]=\FF\cdot sl(V)+\FF'\cdot 1 \ .$$
This proves the first assertion. The second one follows from the obvious fact that the trace $Tr:\FF\cdot
End(V)\to \FF$ is the projection to the second summand of the direct sum decomposition
$$\FF\cdot End(V)=\FF\cdot sl(V)+ \FF '\cdot 1 \ .$$
%This proves the second assertion.
The corollary is proved.
\end{proof}

We can construct more pairs of type $2$ as follows. Let $V$ be a $\kk$-vector space and 
$\Phi:V\times V\to \kk$ be a bilinear form on $V$. Denote by $o(\Phi)$ the orthogonal Lie algebra of $\Phi$, i.e.,
$$o(\Phi)=\{M\in End(V): \Phi(M(u),v)+\Phi(u,M(v))=0~\forall~u,v\in V\} \ .$$
Denote by $K_\Phi\subset V$ the sum of the left and the right
kernels of $\Phi$ (if $\Phi$ is symmetric or skew-symmetric, then
$K_{\Phi}$ is the left kernel of $\Phi$). Finally, denote by
$End(V,K_{\Phi})$ the parabolic subalgebra of $End(V)$ which
consists of all $M\in End(V)$ such that $M(K_{\Phi})\subset
K_{\Phi}$. Clearly, $o(\Phi)\subset End(V,K_{\Phi})$, i.e.,
$(o(\Phi),End(V,K_{\Phi}))$ is a compatible pair. For any object
$\FF$ of $\NN$ we abbreviate
$o(\Phi,\FF):=(o(\Phi),End(V,K_{\Phi}))(\FF)$.

Denote by $sl(V,K_{\Phi})$ the set of all $M$ in $End(V,K_{\Phi})$
such that $Tr(M)=0$ and  $Tr(M_{K_{\Phi}})=0$, where
$M_{K_{\Phi}}:K_{\Phi}\to K_{\Phi}$ is the restriction of $M$ to
$K_{\Phi}$ and $1_{K_{\Phi}}\in End(V,K_{\Phi})$ is any element
such that $\kk\cdot 1+\kk\cdot
1_{K_{\Phi}}+sl(V,K_{\Phi})=End(V,K_{\Phi})$. If $K_{\Phi}=0$, we
set $1_{K_{\Phi}}=0$.

\begin{corollary}
\label{cor:orthogonal and simplectic lie algebra} Let $V\ne 0$ be
a finite-dimensional $\kk$-vector space and $\Phi$ be a symmetric
or skew-symmetric bilinear form on $V$. Then
$(o(\Phi),End(V,K_{\Phi}))$ is of type $2$ and
\begin{equation}
\label{eq:orthogonal and simplectic lie algebra}
o(\Phi,\FF)=\FF\cdot o(\Phi)+\FF'\cdot 1+ \FF'\cdot
1_K+(\FF\FF'+\FF') \cdot sl(V,K_{\Phi})\ .
\end{equation}

\end{corollary}

\begin{proof} We will write $K$ instead of $K_{\Phi}$. First prove that $(\gg,A)=(o(\Phi),End(V,K))$ is of type $2$.
We pass to the algebraic closure of the involved objects, i.e.,
replace both $V$ and $K$ with $\overline V=\overline{\kk}\cdot
V=\overline{\kk}\otimes V$, $\overline K=\overline{\kk}\cdot K$
etc, where $\overline{\kk}$ is the algebraic closure of $\kk$.
Using the obvious fact that $\overline {U+U'}=\overline
U+\overline {U'}$ and $\overline {U\otimes U'}=\overline
{U}\otimes \overline {U'}$ for any subspaces of $End(V)$ and
$\overline{o(\Phi)}=o(\overline{\Phi})$, we see that it suffices
to show that the pair $(o(\overline{\Phi}),End(\overline
V,\overline K)$ is of type $2$.

Furthermore, without loss of generality we consider the case when
$K=0$, i.e., form $\Phi$ is non-degenerate. One can do it by
using the block matrix decomposition with respect to a choice of
compliment of $K$ over which $\Phi$ is nondegenerate.

We will prove the lemma when $\dim V>2$ and leave the rest to the
reader. If $\Phi$ is symmetric, one can choose a basis of
$\overline V$ so that $\overline  V\cong \overline{\kk}^n$, and
$\overline \Phi$ is the standard dot product on
$\overline{\kk}^n$. In this case $o(\overline \Phi)$ is
$o_n(\overline{\kk})$, the Lie algebra of orthogonal matrices,
which is generated by all elements $E_{ij}-E_{ji}$
 where $E_{ij}$ is the corresponding elementary matrix. Using the identity
$(E_{ij}-E_{ji})^2=-(E_{ii}+E_{jj})$ for $i\neq j$, we see that $o_n(\overline{\kk})^2$ contains all diagonal matrices.
Furthermore, if $i,j,k$ are pairwise distinct indices then $(E_{ij}-E_{ji})(E_{jk}-E_{kj})=E_{ik}$.
Thus we have shown that $o_n(\overline{\kk})^2=M_n(\overline{\kk})=End(\overline V,\overline K)$. 
Therefore, $o_n(\kk)^2=M_n(\kk)=End( V, K)$.
This proves the assertion for the symmetric $\Phi$.

If $\Phi $ is skew-symmetric and non-degenerate, then $n=2m$ and one can choose a basis of $\overline V$ such that $V$ is
identified with $\overline{\kk}^n$ and $o(\overline \Phi)$ is identified with the symplectic Lie algebra $sp_{2m}(\overline{\kk})$.
%$\Phi (x,y)=\sum_{i=1}^m x_iy_{m+i}- x_{m+i}y_i$ for $x=(x_1,\ldots,x_{2m})$, $y=(y_1,\ldots,y_{2m})$. Then

Recall that a basis in $sp_{2m}(\overline{\kk})$ can be chosen as
follows. It consists of elements $E_{ij}-E_{j+m,i+m}$,
$E_{i,m+j}+E_{j,m+i}$, $E_{m+i,j}+E_{m+j,i}$, for $i,j\leq m$.
Using the identity $(E_{i,i+m}+E_{i+m,i})^2=E_{ii}+E_{i+m,i+m}$
and the fact that $(E_{ii}-E_{i+m,i+m})\in sp_{2m}(\kk)$, we see
that all diagonal matrices belong
$sp_{2m}(\overline{\kk})+sp_{2m}(\overline{\kk})^2$.

Also, the identity
$(E_{ii}-E_{i+m,i+m})(E_{ij}-E_{j+m,i+m})=E_{ij}$ for $i\neq j$
implies that $E_{ij}\in sp_{2m}(\overline{\kk})^2$ for all
$i,j\leq m$. Similarly, one can prove that $E_{ij}\in
sp_{2m}(\overline{\kk})^2$ for $i,j\geq m$.

Furthermore, the identity $(E_{ii}-E_{i+m,i+m})(E_{i\ell}+E_{i+m,\ell-m})=E_{i\ell}-E_{i+m,\ell-m}$ 
implies that $sp_{2m}(\overline{\kk})+sp_{2m}(\overline{\kk})^2$ contains all $E_{ik}$ for $i\leq m, k>m$ and for
$i>m, k\leq m$. Thus we have shown that 
$sp_{2m}(\overline{\kk})+sp_{2m}(\overline{\kk})^2=M_n(\overline{\kk})=End(\overline V,\overline K)$. 
Therefore, $sp_{2m}(\kk)+sp_{2m}(\kk)^2=M_n(\kk)=End( V, K)$. This proves the assertion for the skew-symmetric $\Phi$.
%Evidently, for all listed Lie algebras $\gg (\kk )^m=gl_n(\kk)$ and $(\gg (\kk )^m)_+=sl_n(\kk )$ for $m\geq 2$. It %follows
%that the second summand in equation \eqref{eq:semisimple center} is equal to $(\FF\FF'+\FF')\cdot sl_n(\kk )$.
%When $\gg (\kk)=sl_n (\kk)$, the second summond is contained in $\FF\cdot sl_n(\kk )$.
%
%Note also that the center of $\langle\gg \rangle$ consists of scalar matrices. It finishes the proof.

Prove \eqref{eq:orthogonal and simplectic lie algebra} now. We
abbreviate $A=End(V,K)$. Recall that $[A,A]=sl(V,K)$ and, if $K\ne
\{0\}$, then $\kk\cdot 1+sl(V,K)$ is of codimension $1$ in $A$,
i.e, $1_K$ always exists. Therefore, applying Proposition
\ref{pr:type two}, we obtain
$$o(\Phi,\FF)=\FF\cdot o(\Phi)+\FF'\cdot A+\FF\FF'[A, A]=\FF\cdot o(\Phi)+\FF'\cdot 1+ \FF'\cdot 1_K+(\FF\FF'+\FF') \cdot sl(V,K)\ .$$
This finishes the proof of Corollary \ref{cor:orthogonal and simplectic lie algebra}.
\end{proof}

Note that our Lie algebras $o(\FF)$ and $sp(\FF)$ do not coincide
with the usual orthogonal and symplectic Lie algebras, which are defined when the  ring $\FF $ posseses an involution.

\section{Upper bounds of $\NN$-current Lie algebras}
\label{sect:upper bounds}
For any compatible pair $(\gg,A)$ define two subspaces $\widetilde{(\gg,A)}(\FF)$ and $\overline{(\gg,A)}(\FF)$ of  $\FF\cdot A$ by:
\begin{equation}
\label{eq:upper bound}\widetilde{(\gg,A)}(\FF)=\FF\cdot \gg+\sum_{k\ge 1} I_k(\FF)
\cdot [\gg,\gg^{k+1}] +[\FF,I_{k-1}(\FF)]\cdot \gg^{k+1} \ ,
\end{equation}
where $I_k(\FF)$ is defined in \eqref{eq:recursion Ikl}; and
\begin{equation}
\label{eq:refined upper bound}
\overline{(\gg,A)}(\FF)=\FF\cdot \gg +
\sum I_{k_1}^{\ell_1+1}I_{k_2}^{\ell_2+1}
\cdot [J_{\ell_1}^{k_1+1},J_{\ell_2}^{k_2+1}]+[I_{k_1}^{\ell_1+1},I_{k_2}^{\ell_2+1}]\cdot J_{\ell_2}^{k_2+1}J_{\ell_1}^{k_1+1},
\end{equation}
where the summation is over all $k_1, k_2\ge 0$, $\ell_1,\ell_2\ge 0$,  and we 
abbreviated $I_k^\ell:=I_k^\ell(\FF)$, $J_k^\ell:=I_k^\ell(A,\gg)$ in notation \eqref{eq:Ikl(S)}.

We will refer to $\widetilde{(\gg,A)}(\FF)$ as the {\it upper bound} of $(\gg,A)(\FF)$ and to
$\overline{(\gg,A)}(\FF)$ as the {\it refined upper bound} of $(\gg,A)(\FF)$.

It is easy to see that the assignments $\FF\mapsto \widetilde{(\gg,A)}(\FF)$ and $\FF\mapsto \overline{(\gg,A)}(\FF)$ 
are functors $\widetilde{(\gg,A)}$ and $\overline{(\gg,A)}$ from $\NN$ to the category $Vect_\kk$ of $\kk$-vector spaces.

The following lemma is obvious.
\begin{lemma} If $(\gg,A)$ is a compatible pair of type $m$ (see Definition \ref{def:type}), then
\begin{equation}
\label{eq:upper bound finite type}
\widetilde{(\gg,A)}(\FF)=\FF\cdot \gg+\sum_{k= 1}^{m-1} I_k(\FF)\cdot [\gg,\gg^{k+1}] +[\FF,I_{k-1}(\FF)]\cdot \gg^{k+1} \
.
\end{equation}

\end{lemma}

Now we formulate the main result of this section, which explains
our terminology and proves that both  $\widetilde{(\gg,A)}$ and
$\overline{(\gg,A)}(\FF)$ define $\NN$-current Lie algebras $\NN\to
{\bf LieAlg}_\NN$.

\begin{theorem}
\label{th:upper bound}
For any compatible pair $(\gg,A)$  and any object $\FF$ of $\NN$  one has:

\noindent (a)
The subspace $\widetilde{(\gg,A)}(\FF)$ is a Lie subalgebra of $\FF\cdot A$.

\noindent (b)
The subspace $\overline{(\gg,A)}(\FF)$ is a Lie subalgebra of $\FF\cdot A$.

\noindent (c) $(\gg,A)(\FF)\subseteq \overline{(\gg,A)}(\FF)\subseteq \widetilde{(\gg,A)}(\FF)$.

\end{theorem}

%%%%%%%%%%%%%%%%%%%%
\begin{proof} We need the following two lemmas.

\begin{lemma}
\label{le:commutators of powers} For any compatible pair
$(\gg,A)$ one has
$$[\gg^{k+1},\gg^m]\subset [\gg,\gg^{k+m}] $$
for any $k,m\ge 1$.
\end{lemma}
This is an obvious consequence of \eqref{eq:Leibniz Jacobi
identity}.

For any subsets $X$ and $Y$  of an object ${\mathcal A}$ of ${\bf
Alg}$ and $\varepsilon\in \{0,1\}$ denote
$$X\bullet_\varepsilon Y:=
\begin{cases} X\cdot Y &\text{if $\varepsilon=0$}\\
 [X, Y] &\text{if $\varepsilon=1$}
\end{cases}$$

\begin{lemma}
\label{le:Gamma grading} Let $\Gamma$ be an abelian group and let
$A$ and $\FF$ be objects of ${\bf Alg}$. Assume that
$E_\alpha\subset \FF$ and $B_\alpha\subset A$ are two families of
subspaces labeled by $\Gamma$ such that
\begin{equation}
\label{eq:Gamma grading} E_\alpha\bullet_\varepsilon
E_\beta\subseteq E_{\alpha+\beta+\varepsilon\cdot v},~B_\beta
\bullet_\varepsilon B_\alpha\subseteq
B_{\alpha+\beta-\varepsilon\cdot v}
\end{equation}
for all $\alpha,\beta\in \Gamma$, $\varepsilon\in \{0,1\}$, where
$v$ is a fixed element of $\Gamma$. Then for any $\alpha_0\in
\Gamma$ the subspace
$$\hh=E_{\alpha_0}\cdot B_{\alpha_0+v}+\sum_{\alpha,\beta\in \Gamma, \varepsilon\in \{0,1\}} (E_\alpha\bullet_{1-\varepsilon} E_\beta)
\cdot (B_{\beta+v}\bullet_\varepsilon B_{\alpha+v})$$ is a Lie
subalgebra of $\FF\cdot A=\FF\otimes A$.
\end{lemma}

\begin{proof}
The equation \eqref{eq:the commutator noncom coeffs} implies that
$$[E\cdot B,E'\cdot B']\subset (E\bullet_{1-\delta} E')\cdot (B'\bullet_\delta B)$$
for each $\delta\in \{0,1\}$.

(i) Set $E=E_\alpha\bullet_{1-\varepsilon} E_\beta$,
$B=B_{\beta+v}\bullet_\varepsilon B_{\alpha+v}$,
$E'=E_{\alpha'}\bullet_{1-\varepsilon'} E_{\beta'}$,
$B'=B_{\beta'+v}\bullet_{\varepsilon'} B_{\alpha'+v}$. Define
$\alpha''=\alpha+\beta+(1-\varepsilon)\cdot  v$ and
$\beta''=\alpha'+\beta'+(1-\varepsilon')\cdot  v$.

Taking into account that $E\subseteq E_{\alpha''}$, $E'\subseteq
E_{\beta''}$, $B\subseteq B_{\alpha''+v}$, and $B'\subseteq
B_{\beta''+v}$ by \eqref{eq:Gamma grading}, we obtain for each
$\delta\in \{0,1\}$:
$$[E\cdot B,E'\cdot B']\subset (E_{\alpha''}\bullet_{1-\delta} E_{\beta''})\cdot (B_{\beta''+v}\bullet_\delta B_{\alpha''+ v})\subset \hh \ .$$

(ii) Set $E=E_{\alpha_0}$, $B=B_{\alpha_0+v}$,
$E'=E_{\alpha'}\bullet_{1-\varepsilon'} E_{\beta'}$,
$B'=B_{\beta'+v}\bullet_{\varepsilon'} B_{\alpha'+v}$. Define
$\beta ''$ as above. Taking into account that $E'\subseteq
E_{\beta''}$ and $B'\subseteq B_{\beta''+v}$ by \eqref{eq:Gamma
grading}, where $\beta''=\alpha'+\beta'+(1-\varepsilon')\cdot v$,
we obtain for each $\delta\in \{0,1\}$:
$$[E\cdot B,E'\cdot B']\subset (E_{\alpha_0}\bullet_{1-\delta} E_{\beta''})\cdot (B_{\beta''+v}\bullet_\delta B_{\alpha_0+ v})\subset \hh \ .$$

(ii) Taking $E=E'=E_{\alpha_0}$, $B=B'=B_{\alpha_0+v}$, we obtain
for each $\delta\in \{0,1\}$:
$$[E\cdot B,E'\cdot B']\subset (E_{\alpha_0}\bullet_{1-\delta} E_{\alpha_0})\cdot (B_{\alpha_0 + v}\bullet_\delta B_{\alpha_0 + v})
\subset \hh \ .$$

The lemma is proved.
\end{proof}

\noindent Now we are going to prove the theorem part-by-part.

 Prove (a).
Using \eqref{eq:the commutator noncom coeffs}, we obtain
$$[\FF\cdot \gg,\FF\cdot \gg]\subset \FF^2\cdot [\gg,\gg]+[\FF,\FF]\cdot \gg^2\subset  \widetilde{(\gg,A)}(\FF)$$
because $\FF^2\subset \FF$,  $[\gg,\gg]\subset \gg$, and $I_0(\FF)=\FF$. Furthermore,
$$[\FF\cdot \gg,I_k(\FF)\cdot [\gg,\gg^{k+1}]]\subset \FF I_k(\FF)\cdot [\gg,[\gg,\gg^{k+1}]]+[\FF,I_k(\FF)]\cdot
[\gg,\gg^{k+1}]\gg\subset  \widetilde{(\gg,A)}(\FF)$$
because $\FF I_k(\FF)\subset I_k(\FF)$,  $[\gg,[\gg,\gg^{k+1}]]\subset [\gg,\gg^{k+1}]$,
and $[\gg,\gg^{k+1}]\gg\subset \gg^{k+2}$. Finally, 
set $J_k:=[\FF,I_{k-1}(\FF)]$.Since $I_{k-1}(\FF)$ is a two-sided ideal in $\FF$, we have
$\FF J_k\subset  \FF I_{k-1}(\FF)\subset I_{k-1}(\FF)$.
Lemma \ref{le:embedded Ikl}(b) taken with $\ell=\infty$, implies $J_k \subset I_k(\FF )$.
Therefore,  $[\FF,J_k]\subset [\FF,I_k(\FF)]$ and 
$$[\FF\cdot \gg,J_k\cdot \gg^{k+1}]\subset \FF\cdot J_k\cdot [\gg,\gg^{k+1}]+
[\FF,J_k]\cdot \gg^{k+2}\subset  \widetilde{(\gg,A)}(\FF).$$

%abbreviating $J_k:=[\FF,I_{k-1}(\FF)]$, we obtain:
%$$[\FF\cdot \gg,J_k\cdot \gg^{k+1}]\subset \FF\cdot J_k\cdot [\gg,\gg^{k+1}]+[\FF,J_k]\cdot \gg^{k+2}\subset  
%\widetilde{(\gg,A)}(\FF)$$
%because, taking into account that $J_k\subset I_{k-1}(\FF)$ 
%$since I_{k-1}(\FF)$ is a two-sided ideal in $\FF$, 
%we have $\FF J_k\subset  \FF I_{k-1}(\FF)\subset I_{k-1}(\FF)$ and, taking
%into account that $J_k\subset I_k(\FF)$ by Lemma \ref{le:embedded
%Ikl}(b) taken with $\ell=\infty$, we have  $[\FF,J_k]\subset $[\FF,I_k(\FF)]$.

Note that for any $k,m\ge 1$ one has:
$$[I_k(\FF)\cdot \gg^{k+1},I_m(\FF)\cdot \gg^{m+1}]\subset $$
$$\subset I_k(\FF)I_m(\FF)\cdot [\gg^{k+1}, \gg^{m+1}]+[I_k(\FF),I_m(\FF)]\cdot \gg^{k+m+2}\subset  \widetilde{(\gg,A)}(\FF)$$
because  $I_k(\FF)I_m(\FF)\subset I_{k+m}(\FF)$ by Lemma
\ref{le:ideals I_k(F) simple}(a), $[\gg^{k+1},\gg^{m+1}]\subset
[\gg,\gg^{k+m}]$ by Lemma \ref{le:commutators of powers}, and
$[I_k(\FF),I_m(\FF)]\subset [\FF,I_{k+m-1}(\FF)]$ by Lemma
\ref{le:ideals I_k(F) simple}(b) taken with $\ell=\infty$.
Therefore, taking into account that
$$[{\bf I}_k,{\bf I}_m]\subset [I_k(\FF)\cdot \gg^{k+1},I_m(\FF)\cdot \gg^{m+1}]\subset \widetilde{(\gg,A)}(\FF)$$
for ${\bf I}_r$ standing for any of the spaces $I_r(\FF)\cdot [\gg,\gg^{r+1}],[\FF,I_{r-1}(\FF)]\cdot \gg^{r+1}$ we finish the proof 
of (a).

%Also denote
%\begin{equation}
%\label{eq:super-upper bound}\overline{\overline{(\gg,A)}(\FF)}=\FF\cdot \gg +\sum_{0\le k_1\le k_2} I_{k_1+k_2}(\FF)\cdot [\gg^{k_1+1},\gg^{k_2+1}]+[I_{k_1}(\FF),I_{k_2}(\FF)]\cdot \gg^{k_1+k_2+2} \ .
%\end{equation}
%
%Clearly, each $I_k(\FF)$ is  a two-sided ideal in $\FF$. It is easy to see that $I_k (\FF)=I_k ^k(\FF)+I_k^{k+1} %(\FF)$.

Prove (b). Taking in Lemma \ref{le:Gamma grading}: $\Gamma=\ZZ^2$,
$\alpha=(k,\ell+1)\in \ZZ^2$, $v=(1,-1)$, one can see that $E_\alpha$
equals to $I_k^{\ell+1}(\FF)$ if $k,\ell\ge 0$ and zero
otherwise. Also, $B_{\alpha+v}$ equals to $I_\ell^{k+1}(A,\gg)$
if $k,\ell\ge 0$ and zero otherwise.

Lemma \ref{le:ideals I_k} implies that   \eqref{eq:Gamma grading}
holds for all $\alpha,\beta\in \ZZ^2$, $\varepsilon\in \{0,1\}$.
Therefore, applying Lemma \ref{le:Gamma grading} with
$\alpha_0=(0,1)$, we  finish the proof of the assertion that
$\overline{(\gg,A)}(\FF)$ is a Lie subalgebra of $\FF\cdot A$.
% This finishes the proof of (b).

Prove (c). The first inclusion $(\gg,A)(\FF)\subset \overline{(\gg,A)}(\FF)$ is 
obvious because  $\FF\cdot \gg\subset \overline{(\gg,A)}(\FF)$ and $\overline{(\gg,A)}(\FF)$ is a Lie subalgebra of $\FF\cdot A$.

Let us prove the second inclusion $(\gg,A)(\FF)\subset \overline{(\gg,A)}(\FF)$ of (c).

Rewrite the result of Lemma \ref{le:ideals I_k(F) simple} (with
$\ell_1=\ell_2=\infty$) in the form of \eqref{eq:Gamma grading}
as:
$$I_{k_1}^{\ell_1+1}(\FF)\bullet_{1-\varepsilon} I_{k_2}^{\ell_2+1}(\FF)\subset I_{k_1}(\FF)\bullet_{1-\varepsilon} I_{k_2}(\FF)
\subset  \begin{cases}
I_{k_1+k_2}(\FF) & \text{if $\varepsilon=1$} \\
[\FF,I_{k_1+k_2-1}(\FF)] & \text{if $\varepsilon=0$}
\end{cases}.$$
Using the obvious inclusion $J_k^{\ell+1}=I_k^{\ell
+1}(A,\gg)\subset \gg^{k+1}$ for all $k,\ell\ge 0$ and Lemma
\ref{le:commutators of powers}, we obtain
$$J_{\ell_2}^{k_2+1}\bullet_{\varepsilon}J_{\ell_1}^{k_1+1}\subset \gg^{k_2+1}\bullet_\varepsilon \gg^{k_1+1}\subset \begin{cases}
\gg^{k_1+k_2+2} & \text{if $\varepsilon=0$} \\
[\gg,\gg^{k_1+k_2+1}] & \text{if $\varepsilon=1$}
\end{cases}
$$
for all $k_1,k_2,\ell_1,\ell_2\ge 0$, $\varepsilon\in \{0,1\}$. Therefore,
we get the inclusion:

$$(I_{k_1}^{\ell_1+1}\bullet_{1-\varepsilon}I_{k_2}^{\ell_2+1})
\cdot (J_{\ell_2}^{k_2+1}\bullet_{\varepsilon}J_{\ell_1}^{k_1+1})\subset $$
$$\subset \begin{cases}
I_{k_1+k_2}(\FF)\cdot [\gg,\gg^{k_1+k_2+1}] & \text{if $\varepsilon=1$} \\
[\FF,I_{k_1+k_2-1}(\FF)]\cdot \gg^{k_1+k_2+2} & \text{if $\varepsilon=0$}
\end{cases} \subset \widetilde{(\gg,A)}(\FF) \ .$$
This proves the inclusion $\overline{(\gg,A)}(\FF)\subset  \widetilde{(\gg,A)}(\FF)$ and finishes the proof of (c).

Therefore, Theorem \ref{th:upper bound} is proved.
\end{proof}

Now we will refine Theorem \ref{th:upper bound} by introducing a
natural filtration on each Lie algebra involved and proving the
``filtered'' version of the theorem.

%%%%%%%%%%%%%%%%%%%%%%

For any compatible pair $(\gg,A)$, any object $\FF$ of $\NN$,
and each $m\ge 1$ we define the subspaces $\FF\cdot \langle
\gg\rangle_m$, $(\gg,A)_m(\FF)$, $\widetilde{(\gg,A)}_m(\FF)$ and
$\overline{(\gg,A)}_m(\FF)$ of $\FF\cdot A$ by:
$$\FF\cdot \langle \gg\rangle_m= \sum_{1\le k\le m} \FF\cdot \gg^{k}$$
\begin{equation}
\label{eq:filtered gg}
(\gg,A)_m(\FF)=\sum_{0\le k< m} (\FF\cdot \gg)^{(k)}
\end{equation}
\begin{equation}
\label{eq:filtered upper bound}
\widetilde{(\gg,A)}_m(\FF)=\FF\cdot \gg+\sum_{1\le k<m} I_k^{\le m-k}(\FF)
\cdot [\gg,\gg^{k+1}] +[\FF,I_{k-1}^{\le m-k}(\FF)]\cdot \gg^{k+1} \ ,
\end{equation}
where $I_k^{\le \ell}(\FF)$ is defined in \eqref{eq:recursion Ikl} and
\begin{equation}
\label{eq:filered refined upper bound}
\overline{(\gg,A)}_m(\FF)=\FF\cdot \gg +
\sum I_{k_1}^{\ell_1+1}I_{k_2}^{\ell_2+1}\cdot [J_{\ell_1}^{k_1+1},J_{\ell_2}^{k_2+1}]+
[I_{k_1}^{\ell_1+1},I_{k_2}^{\ell_2+1}]\cdot J_{\ell_2}^{k_2+1}J_{\ell_1}^{k_1+1},
\end{equation}
where the summation is over all $k_1, k_2\ge 0$, $\ell_1,\ell_2\ge 0$ such that
$k_1+k_2+\ell_1+\ell_2+2\le m$,  and we abbreviated $I_k^\ell:=I_k^\ell(\FF)$, 
$J_k^\ell:=I_k^\ell(A,\gg)$ in the notation \eqref{eq:Ikl(S)}.

Recall that a Lie algebra $\hh=(\hh_1\subset \hh_2\subset \ldots )$ is called a 
{\it filtered Lie algebra} if $[\hh_{k_1},\hh_{k_2}]\subset \hh_{k_1+k_2}$ for all $k_1,k_2\ge 0$.

Taking into account that $[\gg^{k_1+1},\gg^{k_2+1}]\subset
\gg^{k_1+k_2+1}$, one can see that $\hh_m=\FF\cdot \langle
\gg\rangle_m$, $m\ge 1$ defines an increasing filtration on the
Lie algebra $\FF\cdot \langle \gg\rangle$ (where $\langle
\gg\rangle$ is as in \eqref{eq:langle gg rangle}.

The following result is a filtered version of  Theorem \ref{th:upper bound}.
\begin{theorem}
\label{th:filtered upper bound}
For any compatible pair $(\gg,A)$  and an object $\FF$ of $\NN$
\smallskip

\noindent (a)  $\widetilde{(\gg,A)}(\FF)$ is a filtered Lie subalgebra of $\FF\cdot \langle \gg\rangle$.

\noindent (b) $\overline{(\gg,A)}(\FF)$ is a filtered Lie subalgebra of $\FF\cdot \langle \gg\rangle$.

\noindent (c) There is a  chain of inclusions of filtered Lie algebras:
$$(\gg,A)(\FF)\subseteq \overline{(\gg,A)}(\FF)\subseteq \widetilde{(\gg,A)}(\FF)\ .$$

\end{theorem}

The proof of Theorem \ref{th:filtered upper bound} is almost identical to that of Theorem \ref{th:upper bound}.

\section{Perfect pairs and achievable upper bounds}
\label{sect: perfect pairs}

%%%%%%%%%%%%%%%%%%%%%%%%%%%%%%%%%
%%%%%%%%%%%%%%%%%%%%%%%%%%%%%%%%
Below we lay out some sufficient conditions on the compatible pair $(\gg,A)$ which guarantee that the upper bounds are achievable.

\begin{definition}
\label{def:perfect pair} We say that a compatible pair $(\gg,A)$ is {\it perfect} if
\begin{equation}
\label{eq:perfectness}
[\gg,\gg^k]\gg+(\gg^k\cap \gg^{k+1})=\gg^{k+1}
\end{equation}
for all $k\ge 2$.

%We say that a compatible pair $(\gg,A)$ is split semisimple if $\gg$ is a split semisimple (over $\kk$) Lie subalgebra of $A$.

\end{definition}

\begin{definition}
\label{def:strongly graded}
We say that a  Lie algebra $\overline \gg$ over an algebraically closed field $\overline{\kk}$ is {\it strongly graded} if there exists an element $h_0\in \overline{\gg}$ such that:

\noindent (i) The operator $ad~h_0$ on $\overline{\gg}$ is diagonalizable, i.e.,

\begin{equation}
\label{eq:h-decomposition}
\overline{\gg}=\bigoplus_{c\in \overline{\kk}} \overline{\gg}_c \ ,
\end{equation}
where  $\overline{\gg}_c\subset \overline \gg$ is an eigenspace of $ad~h_0$ with the eigenvalue $c$.

\noindent (ii) The nulspace $\overline{\gg}_0$ of $ad~h_0$ is spanned by  $[\overline{\gg}_c,\overline{\gg}_{-c}]$, $c\in \overline{\kk}\setminus\{0\}$.

\end{definition}

The class of strongly graded Lie algebras is rather large; it
includes all semisimple and Kac-Moody Lie algebras, as well as
the Virasoro algebra.
%\comment{$W$-algebras}

\begin{maintheorem}
\label{th:perfect} Let $(\gg,A)$ be a compatible pair. Then

\noindent (a)  If $(\gg,A)$ is perfect, then for any object $\FF$ of $\NN$ one has
$$(\gg,A)(\FF)=\widetilde{(\gg,A)}(\FF) \ ,$$
i.e., the $\NN$-current Lie algebras $(\gg,A),\widetilde{(\gg,A)}:\NN\to {\bf LieAlg}_\NN$ are equal.

\noindent (b) If $\overline \gg=\overline{\kk}\otimes \gg$  is strongly graded, then $(\gg,A)$ is perfect.

\noindent (c)  If $\gg$ is  semisimple over $\kk$, then for any object $\FF$ of $\NN$ one has
\begin{equation}
\label{eq:semisimple center}
(\gg,A)(\FF)=\FF\cdot \gg+\sum_{k\ge 2} I_{k-1}(\FF)\cdot  (\gg^k)_+ + [\FF,I_{k-2}(\FF)]\cdot Z_k(\gg)\ ,
\end{equation}
where $(\gg^k)_+=[\gg,\gg^k]$ is
the ``centerless'' part of $\gg^k$, $Z_k(\gg)=Z(\langle\gg\rangle)\cap
\gg^k$, and $Z(\langle\gg\rangle)$ is the center of
$\langle\gg\rangle=\sum_{k\ge 1} \gg^k$.

\end{maintheorem}

\begin{proof} To prove the theorem we need a proposition and two lemmas.

\begin{proposition}
\label{pr: propagation of gg} Let $(\gg,A)$ be a compatible pair
and $\FF$ be an object of $\NN$.

\noindent (a) Assume that for some $k\ge 1$ one has
$$ I\cdot [\gg,\gg^k] \subset (\gg,A)(\FF) $$
where  $I$ is a left ideal in $\FF$. Then:
\begin{equation}
\label{eq:kth commutator belongs to gg} [\FF,I]\cdot [\gg,\gg^k]
\gg\subset (\gg,A)(\FF) \ .
\end{equation}

\noindent (b) Assume that for some $k\ge 1$ one has
$$ J\cdot \gg^k \subset (\gg,A)(\FF) $$
where  $J$ is a subset of $\FF$ such that  $[\FF,J]\subset J$.
Then:

%\noindent (a)  For any  one has
%\begin{equation}
%\label{eq:kth commutator belongs to gg}
%\FF^{(k)}\cdot \gg^{k+1}\subset (\gg,A)(\FF)
%\end{equation}
%for all $k\ge 0$;
%\begin{equation}
%\label{eq:k1th times k2nd commutator belong to gg}
%\FF^{(k_1)}\FF^{(k_2)}\cdot [\gg^{k_1+1},\gg^{k_2+1}]\subset (\gg,A)(\FF)
%\end{equation}
%for all $k_1,k_2\ge 0$.
\begin{equation}
\label{eq:kth power belongs to gg} [\FF,J]\cdot \gg^{k+1}+(\FF
J+J)\cdot [\gg,\gg^k]\subset (\gg,A)(\FF)
\end{equation}

\end{proposition}

\begin{proof} Prove (a). Indeed,
$$[\FF\cdot \gg,I\cdot [\gg,\gg^k]]\equiv  [\FF,I]\cdot [\gg,\gg^k]\gg \mod \FF I\cdot [\gg,[\gg,\gg^k]] \ .$$
Since $\FF I\subset I$ and $[\gg,[\gg,\gg^k]]\subset [\gg,\gg^k]$,
and, therefore, $\FF I\cdot [\gg,[\gg,\gg^k]]\subset I\cdot
[\gg,\gg^k] \subset (\gg,A)(\FF)$, the above congruence implies
that $[\FF,I]\cdot [\gg,\gg^k]\gg$ also belongs to
$(\gg,A)(\FF)$. This proves (a).

Prove (b).  For any $g\in \gg$ we obtain:
$$[\FF\cdot g,J\cdot g^k]=[\FF,J]\cdot g^{k+1}$$
which implies that $[\FF,J]\cdot \widetilde {\gg^{k+1}}\subset
(\gg,A)(\FF)$ (in the notation of Lemma \ref{le:powers of gg}).
Using Lemma \ref{le:powers of gg}, we obtain
$$[\FF,J]\cdot \widetilde {\gg^{k+1}}\equiv [\FF,J]\cdot \gg^{k+1}\mod [\FF,J]\cdot (\gg^k\cap \gg^{k+1}) \ .$$
Taking into account that $[\FF,J]\cdot (\gg^k\cap
\gg^{k+1})\subset [\FF,J]\cdot \gg^k \subset (\gg,A)(\FF)$, the
above formula implies that $[\FF,J]\cdot \gg^{k+1}$ also belongs to
$(\gg,A)(\FF)$. Furthermore,
$$[\FF\cdot \gg,J\cdot \gg^k]\equiv \FF J\cdot [\gg,\gg^k] \mod [\FF,J]\cdot \gg^{k+1} \ .$$
Therefore, using the already proved inclusion $[\FF,J]\cdot
\gg^{k+1}\subset (\gg,A)(\FF)$, we see that $\FF J\cdot
[\gg,\gg^k]$ also belongs to $(\gg,A)(\FF)$. Finally, using the
fact that $[\gg,\gg^k]\subset \gg^k$, we obtain $J\cdot
[\gg,\gg^k]\subset J\cdot \gg^k\subset (\gg,A)(\FF)$. This proves
(b).

Proposition \ref{pr: propagation of gg} is proved.
\end{proof}

\begin{lemma}
\label{le:h-decomposition plus cartan} Let $(\gg,A)$ be a
compatible pair. Assume that $h_0\in \overline
\gg=\overline{\kk}\otimes \gg$ is such that  $ad~h_0$ is
diagonalizable, i.e., has a decomposition
\eqref{eq:h-decomposition}. Then:

\noindent (a)  For each $k\ge 1$ and each  ${\bf
c}=(c_1,\ldots,c_{k+1})\in \overline{\kk}^{k+1}\setminus \{0\}$
the subspace $\overline \gg_{c_1}\cdots \overline \gg_{c_{k+1}}$
of $\overline \gg^{k+1}$ belongs to $[\overline \gg,\overline
\gg^k]\overline \gg+(\overline \gg^k\cap \overline \gg^{k+1})$.

\noindent (b) If $\overline \gg=\overline{\kk}\otimes \gg$ is a
strongly graded Lie algebra, then  one has (in the notation of
Definition \ref{def:strongly graded}):
$$\overline \gg_0^{k+1}\subset [\overline \gg,\overline \gg^k]\overline \gg+(\overline \gg^k\cap \overline \gg^{k+1}) \ .$$

\end{lemma}

\begin{proof} Prove (a). Clearly, under the adjoint action of $h_0$ on $\overline \gg^k$ each vector of
$x\in \overline \gg_{c_1}\cdots \overline \gg_{c_k}$ satisfies
$[h_0,x]=(c_1+\cdots + c_k)x$. Therefore, for any
$(c_1,\ldots,c_k)\in \overline{\kk}^k$ such that $c_1+\cdots
+c_k\ne 0$ the subspace $\overline \gg_{c_1}\cdots \overline
\gg_{c_k}$ belongs to $[\overline \gg,\overline \gg^k]$. Clearly,
$$\overline \gg_{c_1}\cdots \overline \gg_{c_{k+1}}\equiv \overline \gg_{c_{\sigma(1)}}\cdots \overline \gg_{c_{\sigma(k)}}
\overline \gg_{c_{\sigma(k+1)}}\mod \overline \gg^k\cap \overline
\gg^{k+1} $$ for any  permutation $\sigma\in S_{k+1}$.

It is also easy to see that for any  ${\bf
c}=(c_1,\ldots,c_{k+1})\in \overline{\kk}^{k+1}\setminus \{0\}$
there exists a permutation $\sigma\in S_{k+1}$ such that
$c_{\sigma(1)}+\cdots +c_{\sigma(k)}\ne 0$ and, therefore,
$$\overline \gg_{c_1}\cdots \overline \gg_{c_{k+1}}\subset (\overline \gg_{c_{\sigma(1)}}\cdots
\overline \gg_{c_{\sigma(k)}})\overline \gg_{c_{\sigma(k+1)}}+
\overline \gg^k\cap \overline \gg^{k+1}\in [\overline
\gg,\overline \gg^k]\overline \gg + (\overline \gg^k\cap
\overline \gg^{k+1}) \ .$$
 This proves (a).

 Prove (b) now.  There exists $c\in \overline{\kk}\setminus \{0\}$ such that $\overline \hh_c=[\overline \gg_c,\overline \gg_{-c}]\ne 0$.
% choose $E\in \gg_c$, $F\in \gg_c$ such that $h_c=[E,F]\ne 0$.
Then  we obtain the following congruence:
$$[\overline \gg_c,\overline \gg_0^{k-1}\overline \gg_{-c}]\gg_0\equiv \overline \gg_0^{k-1}
\overline \hh_c\overline \gg_0 \mod [\overline \gg_c,\overline
\gg_0^{k-1}]\overline \gg_{-c}\overline \gg_0$$ Taking into
account that $[\overline \gg_c,\overline \gg_0]\subset \overline
\gg_c$, we obtain:
$$[\overline \gg_c,\overline \gg_0^{k-1}]\overline \gg_{-c}\overline \gg_0\subset  \sum_{i=1}^{k-1}
\overline \gg_0^{i-1}\overline \gg_c \overline
\gg_0^{k-1-i}\overline \gg_{-c} \overline \gg_0\subset [\overline
\gg,\overline \gg^k]\overline \gg + (\overline \gg^k\cap
\overline \gg^{k+1})$$ by the already proved part (a). Therefore,
$\overline \gg_0^{k-1} \overline \hh_c\overline \gg_0\subset
[\overline \gg,\overline \gg^k] \overline \gg + (\overline
\gg^k\cap \overline \gg^{k+1})$. Since $\overline \gg$ is
strongly graded, the subspaces $\hh_c$, $c\in
\overline{\kk}\setminus\{0\}$ span $\overline \gg_0$, and
therefore,  $\overline \gg_0^{k+1}\subset [\overline
\gg,\overline \gg^k]\overline \gg + (\overline \gg^k\cap
\overline \gg^{k+1})$. This proves (b).

The lemma is proved.
\end{proof}

\begin{lemma}
\label{le:semisimple enveloping center} Let $\gg$ be a semisimple
Lie algebra over $\kk$. Then for any compatible pair $(\gg, A)$
one has the following decomposition of the $\gg$-module $\gg^k$,
$k\ge 2$:
$$\gg^k=[\gg,\gg^k]+ Z_k(\gg),~[\gg,\gg^k]\cap  Z_k(\gg)=\{0\} \ ,$$
where $Z_k(\gg)=Z(\langle\gg\rangle)\cap \gg^k$, and
$Z(\langle\gg\rangle)$ is the center of
$\langle\gg\rangle=\sum_{k\ge 0} \gg^k$.
\end{lemma}

\begin{proof} Clearly, $\gg^k$ is a semisimple finite-dimensional $\gg$-module (under the adjoint action).
Therefore, it uniquely decomposes into isotypic components one of
which, the component of invariants, is $Z_k(\gg)$. Denote the sum
of all non-invariant isotypic components by $(\gg^k)_+$.  By
definition, $\gg^k=(\gg^k)_++Z_k(\gg)$ and $(\gg^k)_+\cap
Z_k(\gg)=0$. It remains to prove that $(\gg^k)_+=[\gg,\gg^k]$.
Indeed, $[\gg,\gg^k]\subseteq (\gg^k)_+$. On the other hand, each
non-trivial irreducible $\gg$-submodule $V\subset \gg^k$ is
faithful, i.e., $[\gg,V]=V$ (since $[\gg,V]$ is always a
$\gg$-submodule of $V$). Therefore, $[\gg,\gg^k]$ contains all
non-invariant isotypic components, i.e., $[\gg,\gg^k]\subset
(\gg^k)_+$. The double inclusion obtained implies that
$(\gg^k)_+=[\gg,\gg^k]$. The lemma is proved.
\end{proof}

\noindent Now we are ready to prove the theorem part-by-part.

Prove (a).
In view of Theorem \ref{th:upper bound}(c), it suffices to prove
that $\widetilde{(\gg,A)}(\FF)\subset (\gg,A)(\FF)$, that is,
\begin{equation}
\label{eq:lower perfect}
I_k(\FF)\cdot [\gg,\gg^{k+1}]\subset (\gg,A)(\FF),~[\FF,I_k(\FF)]\cdot\gg^{k+2}\subset (\gg,A)(\FF)
\end{equation}
for $k\ge 0$.

We will prove \eqref{eq:lower perfect} by induction on $k$. First,
verify the base of induction at $k=0$. Obviously, $I_0(\FF)\cdot
[\gg,\gg]\subset FF\cdot [\gg,\gg]\subset (\gg,A)(\FF)$.
Furthermore, Proposition \ref{pr: propagation of gg}(b) taken
with $k=1$, $J=\FF$ implies that $[\FF,\FF]\cdot \gg^2\subset
(\gg,A)(\FF)$.

Now assume that $k>0$. Using a part of the inductive hypothesis
in the form $[\FF,I_{k-1}(\FF)]\cdot\gg^{k+1}\subset
(\gg,A)(\FF)$ and applying Proposition \ref{pr: propagation of
gg}(b) with $J=[\FF,I_{k-1}(\FF)]$, we obtain $(\FF J+J)\cdot
[\gg,\gg^{k+1}]\subset (\gg,A)(\FF)$. In its turn, Lemma
\ref{le:embedded Ikl}(c) taken with $\ell=\infty$ implies that
$\FF J+J=I_k(\FF)$. Therefore, we obtain
$$I_k(\FF)\cdot [\gg,\gg^{k+1}]\subset (\gg,A)(\FF)\ ,$$
which is the first inclusion of \eqref{eq:lower perfect}. To prove the second inclusion \eqref{eq:lower perfect}, we will use Proposition \ref{pr: propagation of gg}(a) with $I=I_k(\FF)$:
$$(\gg,A)(\FF)\supset [\FF,I_k(\FF)]\cdot [\gg,\gg^{k+1}]\gg \ .$$
On the other hand, using the perfectness of the pair $(\gg,A)$, we obtain:
$$[\FF,I_k(\FF)]\cdot [\gg,\gg^{k+1}]\gg\equiv[\FF,I_k(\FF)]\cdot \gg^{k+2}\mod [\FF,I_k(\FF)]\cdot \gg^{k+1}\cup \gg ^{k+2} \ .$$
But Lemma \ref{le:embedded Ikl}(a) taken with $\ell=\infty$ implies that $I_k(\FF)\subset I_{k-1}(\FF)$, therefore,
$$[\FF,I_k(\FF)]\cdot(\gg^{k+1}\cap \gg^{k+2})\subset [\FF,I_k(\FF)]\cdot\gg^{k+1} \subset [\FF,I_{k-1}(\FF)]\cdot\gg^{k+1}\subset (\gg,A)(\FF)$$
by the inductive hypothesis. This gives the second inclusion of \eqref{eq:lower perfect}. Therefore, Theorem \ref{th:perfect}(a) is proved.

Prove (b) now.
%We will first show that each the pair $(\gg,A)$ is perfect, whenever $\overline \gg=\overline{\kk}\otimes \gg$ is strongly graded.
%We need the following result.
Lemma \ref{le:h-decomposition plus cartan} guarantees that for
any strongly graded Lie algebra $\gg$ one has:
$$[\overline \gg,\overline \gg^k]\overline \gg+(\overline \gg^k\cap \overline \gg^{k+1})=\overline \gg^{k+1}$$
for all $k\ge 2$, where $\overline \gg=\overline{\kk}\otimes \gg$
is the ``algebraic closure'' of $\gg$. Since the ``algebraic
closure'' commutes with the multiplication and the commutator
bracket in $A$, the restriction of the above equation to
$\gg^{k+1}\subset \overline \gg^{k+1}$ becomes
\eqref{eq:perfectness}.

This finishes the proof of Theorem \ref{th:perfect}(b).

\noindent Prove (c) now. Since for each semisimple Lie algebra $\gg$ 
the compatible pair $(\gg,A)$ is perfect by the already proved 
Theorem \ref{th:perfect}(b),  Theorem \ref{th:perfect}(a) 
implies that $(\gg,A)(\FF)=\widetilde {(\gg,A)}(\FF)$. Therefore, in order to finish 
the proof of Theorem \ref{th:perfect}(c), it suffices to show that
\begin{equation}
\label{eq:semisimple center tilde}
\widetilde {(\gg,A)}(\FF)=\FF\cdot \gg+\sum_{k\ge 2} I_{k-1}(\FF)\cdot  (\gg^k)_+ + [\FF,I_{k-2}(\FF)]\cdot Z_k(\gg) \ .
\end{equation}
%We need the following simple fact regarding compatible pairs $(\gg,A)$ with $\gg$ semisimple.

Using Lemma \ref{le:semisimple enveloping center} and the
definition \eqref{eq:upper bound} of $\widetilde {(\gg,A)}(\FF)$,
we obtain
$$\widetilde{(\gg,A)}(\FF)=\FF\cdot \gg+\sum_{k\ge 1} I_k(\FF)\cdot [\gg,\gg^{k+1}] +[\FF,I_{k-1}(\FF)]\cdot \gg^{k+1}$$
$$= \FF\cdot \gg+\sum_{k\ge 1} (I_k(\FF)+[\FF,I_{k-1}(\FF)])\cdot [\gg,\gg^{k+1}] +[\FF,I_{k-1}(\FF)]\cdot Z_{k+1}(\gg)\ ,$$
which, after taking into account that $[\FF,I_{k-1}(\FF)]\subset I_k(\FF)$ 
(and shifting the index of summation), becomes the right hand side of 
\eqref{eq:semisimple center tilde}. This finishes the proof of Theorem \ref{th:perfect}(b).

Therefore, Theorem \ref{th:perfect} is proved.
\end{proof}

The following is a direct corollary of Theorem \ref{th:perfect}.

\begin{corollary} Assume that a compatible pair $(\gg,A)$ is perfect and  $\FF$ is a $\kk$-algebra satisfying $I_1(\FF)=\FF$. Then
\begin{equation}
\label{eq:simple F}
(\gg,A)(\FF)= \FF\cdot \gg+\FF\cdot [\gg,\langle \gg\rangle]+[\FF,\FF]\cdot \langle \gg\rangle
\end{equation}
(where $\langle \gg\rangle=\sum_{k\ge 1} \gg^k$).
\end{corollary}

\begin{proof} First, show by induction that $I_k(\FF)=\FF$ for all $k\ge 1$.
It follows immediately from Lemma \ref{le:embedded Ikl}(c)
implying that $I_{k+1}(\FF)=\FF[\FF, I_k(\FF)] + [\FF, I_k(\FF)]$.

%If $k=1$, we have nothing to prove. Furthermore, using the
%inclusion $I_{k-1}(\FF)\FF'\subseteq I_k(\FF)$ and the inductive
%hypothesis $I_{k-1}(\FF)=\FF$, we obtain
%$\FF=\FF\FF'=I_{k-1}(\FF)\FF'\subseteq I_k(\FF)$, therefore,
%$I_k(\FF)=\FF$.
This and \eqref{eq:upper bound} imply that
$$\widetilde {(\gg,A)}(\FF)= \FF\cdot \gg +\sum_{k\ge 1} \FF\cdot [\gg,\gg^{k+1}] +[\FF,\FF]\cdot \gg^{k+1}=
\FF\cdot \gg+\FF\cdot [\gg,\langle \gg\rangle]+[\FF,\FF]\cdot \langle \gg\rangle  \ .$$
This and Theorem \ref{th:perfect} finish the proof.
\end{proof}

\begin{remark} The condition $I_1(\FF)=\FF$ holds for each noncommutative simple unital algebra $\FF$, 
e.g., for each noncommutative skew-field $\FF$ containing $\kk$. 
Therefore, for all such algebras and any perfect pair $(\gg,A)$, 
the Lie algebra $(\gg,A)(\FF)$ is given by the relatively simple formula 
\eqref{eq:simple F}, which also complements \eqref{eq:type two}.

\end{remark}

The following result is a specialization of Theorem \ref{th:perfect} to the case when $\gg=sl_2(\kk)$.
\begin{theorem}
\label{th:perfect sl_2} Let $A$ an object of ${\bf Alg}_1$ containing $sl_2(\kk)$ as a Lie subalgebra. Then
\begin{equation}
\label{eq:semisimple center sl_2}
(sl_2(\kk),A)(\FF)=\FF\cdot sl_2(\kk)+[\FF\cdot 1,Z_1(A,\FF)] +\sum_{i\ge 1} Z_i(A,\FF)\cdot  V_{2i}\ ,
\end{equation}
where
$$Z_i(A,\FF)=\sum_{j\ge 0} I_{i+2j-1}(\FF)\cdot\Delta^j\ ,$$
$\Delta=2EF+2FE+H^2$ is the Casimir element, and $V_{2i}$ is the $sl_2(\kk)$-submodule of $A$ generated by $E^i$.
In particular, if $A=End(V)$, where $V$ is a simple $(m+1)$-dimensional $sl_2(\kk)$-module, then
\begin{equation}
\label{eq:semisimple center sl_2 fd}
(sl_2(\kk),A)(\FF)=[\FF,\FF]\cdot 1+\sum_{k=1}^m I_{k-1}(\FF)\cdot  V_{2k} \ .
\end{equation}

\end{theorem}

\begin{proof} Prove \eqref{eq:semisimple center sl_2}. Clearly, each $\gg^k$ is a finite-dimensional $sl_2(\kk)$-module generated by the highest weight vectors $\Delta^jE^i$, $i,j\ge 0$ , $i+2j\le k$. That is, in notation of \eqref{eq:semisimple center}, one has
$$(\gg^k)_+=\sum_{i>0,j\ge 0 , i+2j\le k} \Delta^j\cdot V_{2i} \ , $$
where the sum is direct (but some summands may be zero) and
$$V_{2i}=\sum_{r=-i}^i \kk\cdot (ad~F)^{i+r}(E^i)$$ is the corresponding simple $sl_2(\kk)$-module; and
$$Z_k(\gg)=\sum_{1\le j\le k/2} \kk\cdot \Delta^j\ ,$$
where the sum is direct. Therefore, taking into account that $I_k(\FF)\subset I_{k-1}(\FF)$, the equation
\eqref{eq:semisimple center} simplifies to
$$(sl_2(\kk),A)(\FF)=\FF\cdot \gg+\sum_{i>0,j\ge 0} I_{i+2j-1}(\FF)\cdot  \Delta^j V_{2i} + 
\sum_{j\ge 1}[\FF,I_{2j-2}(\FF)]\cdot \Delta^j$$
$$=\FF\cdot sl_2(\kk)+[\FF\cdot 1,Z_1(A,\FF)] +\sum_{i\ge 1} Z_i(A,\FF)\cdot  V_{2i}\ .$$
This finishes the proof of \eqref{eq:semisimple center sl_2}.

Prove \eqref{eq:semisimple center sl_2 fd}. Indeed, now  $\Delta\in \kk\setminus \{0\}$, $E^k=0$ for $k>m$. Therefore
$$Z_k(A,\FF)=\sum_{j\ge 0} I_{k+2j-1}(\FF)\cdot\Delta^j=\sum_{j\ge 0} I_{k+2j-1}(\FF)\cdot 1=I_{k-1}(\FF)\cdot 1$$
because $I_k(\FF)\subset
I_{k-1}(\FF)$. Finally, using  \eqref{eq:semisimple center sl_2}, we obtain:
$$(sl_2(\kk),A)(\FF)=\FF\cdot sl_2(\kk)+[\FF,Z_1(A,\FF)]\cdot 1 +\sum_{i\ge 1} Z_i(A,\FF)\cdot  V_{2i}$$
$$=[\FF,\FF]\cdot 1+\sum_{1\le k\le m} I_{k-1}(\FF)\cdot  V_{2k} \ .$$
Theorem \ref{th:perfect sl_2} is proved.
\end{proof}

\bigskip

\section{$\NN$-groups}
\label{sect: n-groups}

Throughout the section we assume that each object of $\NN$ is  a unital $\kk$-algebra, i.e.,  $\NN$ is a sub-category of ${\bf Alg}_1$.

\subsection{From $\NN$-Lie algebras to $\NN$-groups and generalized $K_1$-theories}
\label{From NN-Lie algebras to NN-groups and generalized K1-theories}
In this section we use $\FF$-algebras and the category ${\bf LieAlg}_\NN$ defined in Section \ref{sect:Noncommutative current Lie
algebras}.
\begin{definition}
An {\it affine $\NN$-group} is a triple $(\FF,\GG,{\mathcal A})$,
where $\FF$ is an object of $\NN$, ${\mathcal A}$ is an
$\FF$-algebra in ${\bf Alg}_1$ (i.e., $\iota:\FF\to {\mathcal A}$
respects the unit), and $\GG$ is a subgroup of the group of units
${\mathcal A}^\times$ such that $\GG$ contains the image
$\iota(\FF^\times)=\iota(\FF)^\times$. A morphism
$(\FF_1,\GG_1,{\mathcal A}_1)\to (\FF_2,\GG_2,{\mathcal A}_2)$ of
affine $\NN$-groups is  a pair $(\varphi,\psi)$, where
$\varphi:\FF_1\to \FF_2$ is a  morphism in $\NN$ and
$\psi:\GG_1\to \GG_2$ is a group homomorphism such that
$\psi\circ \iota_1|_{\FF_1^\times}=\iota_2|_{\FF_1^\times}\circ
\psi$.

Denote by ${\bf Gr}_\NN$ the category of affine $\NN$-groups.

%Note that if $\NN=(\FF,Id_\FF)$ has only one object $\FF$ and only the identity arrow $Id_\FF$, then an $\NN$-group is simply any object of ${\bf %Gr}_\NN$ of the form $(\FF,\GG,{\mathcal A})$.

%Let $\FF$ be an object of ${\bf Alg}_1$, we say that an algebra in the category $\FF-Mod-\FF$ in is {\it unital} algebra  if one has the unit  %$1_\FF:\FF\hookrightarrow {\mathcal A}$ such that the $\FF$-actions become multiplications with the copy of $\FF$ in ${\mathcal A}$.

Next, we will construct a number of affine $\NN$-groups out of a given affine $\NN$-group or a given $\NN$-algebra as follows.

Let ${\bf LieAlg}_{\NN;1}$ be the sub-category of ${\bf LieAlg}_\NN$ whose objects are triples $(\FF,{\mathcal L},{\mathcal A})$, where $\FF$ is an object of ${\bf Alg}_1$, ${\mathcal A}$ is an $\FF$-algebra in ${\bf Alg}_1$, and ${\mathcal L}$ is a Lie subalgebra of ${\mathcal A}$ invariant under the adjoint action of $\iota(\FF)$; morphisms in ${\bf LieAlg}_{\NN;1}$ are those morphisms in ${\bf LieAlg}_\NN$ which respect the unit.

Given an object $(\FF,{\mathcal L},{\mathcal A})$  of ${\bf LieAlg}_{\NN;1}$, define the triple
$Exp(\FF,{\mathcal L},{\mathcal A}):=(\FF,\GG,{\mathcal A})$, where $\GG$ is the subgroup of ${\mathcal A}^\times$ generated by $\iota(\FF)^\times$ and by the stabilizer $\{g\in {\mathcal A}^\times:g{\mathcal L}g^{-1}={\mathcal L}\}$ of ${\mathcal L}$ in  ${\mathcal A}^\times$.

Given an object $(\FF,\GG,{\mathcal A})$  of ${\bf Gr}_\NN$, define the triple
$Lie(\FF,\GG,{\mathcal A}):=(\FF,{\mathcal L},{\mathcal A})$, where ${\mathcal L}$ is the Lie subalgebra of ${\mathcal A}$ generated (over $\kk$) by the set $\{g\cdot \iota(f)\cdot g^{-1}:g\in \GG, f\in \FF\}$, that is,  ${\mathcal L}$ is the smallest Lie subalgebra of ${\mathcal A}$  containing $\iota(\FF)$ and invariant under conjugation by $\GG$ (therefore, ${\mathcal L}$ is invariant under the adjoint action of the subalgebra $\iota(\FF)$).

The following result is obvious.
\begin{lemma}
\label{le:N-Lie correspondence}
For each object  $(\FF,{\mathcal L},{\mathcal A})$ of ${\bf LieAlg}_{\NN;1}$ the triple $Exp(\FF,{\mathcal L},{\mathcal A})$ is an affine $\NN$-group; and 
for any affine $\NN$-group $(\FF,\GG,{\mathcal A})$ the triple $Lie(\FF,\GG,{\mathcal A})$ is an object of ${\bf LieAlg}_{\NN;1}$.

\end{lemma}

\begin{remark} The operations $Exp$ and $Lie$ are analogues of the Lie correspondence (between Lie algebras and Lie groups). However, similarly to the operations $L_i$ from from Lemma \ref{le:Lie functors}, in general they are not functors.
\end{remark}

Composing these two operations with each other and the operations
$L_i$, we can obtain a 
number of affine $\NN$-groups out of a given $\NN$-Lie algebra and
vice versa.

%$\bullet$  For each object $(\FF,\GG,{\mathcal A})$ of ${\bf Gr}_\NN$ denote  $G_{max}(\FF,\GG,{\mathcal A}):=(\FF,\GG_{max},{\mathcal A})=Exp(Lie(\FF,\GG,{\mathcal A}))$. By the construction, %$\GG\subset \GG_{max}$.
%
%
%$\bullet$ For each object $(\FF,{\mathcal L},{\mathcal A})$ of ${\bf LieAlg}_\NN^1$ denote  $L_{min}(\FF,{\mathcal L},{\mathcal A}):=(\FF,{\mathcal L}_{min},{\mathcal A})$, where ${\mathcal %L}_{min}={\mathcal L}\cap {\mathcal L}'$, where $(\FF,{\mathcal L}',{\mathcal A})=Lie(Exp(\FF,{\mathcal L},{\mathcal A}))$. By the construction,  ${\mathcal L}'$ is an $\FF$-Lie subalgebra of the %$\FF$-Lie algebra ${\mathcal L}+\iota(\FF)$; and ${\mathcal L}_{min}$ is an $\FF$-Lie subalgebra of the $\FF$-Lie algebra ${\mathcal L}$.
%
%\begin{proposition} The above correspondences define functors $L_{min}:{\bf LieAlg}_\NN^1\to {\bf LieAlg}_\NN^1$ and $G_{max}:{\bf Gr}_\NN\to {\bf Gr}_\NN$  and ``injective''
%natural transformations $L_{min}\to Id$ and $Id\to G_{max}$ respectively.
%
%\end{proposition}
%
%
%We can  construct new  $\NN$-Lie groups and new $\NN$-Lie algebras by composing the above functors.
%
%\begin{remark} The functor $L_{min}$ is a ``purely noncommutative'' construction, because for any object  $(\FF,{\mathcal L},{\mathcal A})$ of ${\bf LieAlg}_\NN^1$ with a commutative algebra $\FF$, one %obtains $L_{min}(\FF,{\mathcal L},{\mathcal A})=(\FF,\{0\},{\mathcal A})$, i.e., the  trivial object.
%\end{remark}

By definition, one has a natural (forgetful) projection functor $\pi:{\bf
Gr}_\NN\to \NN$ by $\pi(\FF,\GG,{\mathcal A})=\FF$ and $\pi(\varphi,\psi)=\varphi|_{\FF_1}$.

A {\it noncommutative current group} (or simply {\it $\NN$-current group}) is any functor ${\mathfrak G}:\NN\to {\bf Gr}_\NN$ such that
$\pi\circ {\mathfrak G}=Id_\NN$ (i.e., ${\mathfrak G}$ is a section of $\pi$).
\end{definition}

Note that if $\NN=(\FF,Id_\FF)$ has only one object $\FF$ and only the identity arrow $Id_\FF$, 
then the $\NN$-current group is simply any object of ${\bf Gr}_{{\bf Alg}_1}$ of the
form
$(\FF,{\mathcal G},{\mathcal A})$. In this case, we will sometimes refer to $\GG$ an {\it $\FF$-current group}.

The above arguments allow for constructing a number of $\FF$-current groups out of $\FF$-current 
Lie algebras and vice versa. In a a different situation, for any subcategory $\NN$ of ${\bf Alg}_1$, we will
construct below
a class of $\NN$-current groups associated with compatible pairs $(\gg,A)$. 
More general $\NN$-current groups will be considered elsewhere.

Similarly to Section \ref{sect:Noncommutative current Lie algebras}, given an object $\FF$ of ${\bf Alg}_1$ and a group
$G$, we refer to a group homomorphism $\iota:\FF^\times\to G$ in  ${\bf Alg}$ as a  {\it $\FF$-group structure on $G$} (we will also 
refer to $G$ an {\it  $\FF$-group}).

%%%%%%%%%%%%%%%%%%%
\begin{definition} A {\it decorated group} is a pair $(\FF,\GG)$, where $\FF$ is an object of ${\bf Alg}_1$ and $\GG$ is an $\FF$-group.

We denote by ${\bf DecGr}$ the category whose objects are decorated groups and morphisms are pairs $(\varphi,\psi):(\FF_1,\GG_1)\to (\FF_2,\GG_2)$,
where $\varphi:\FF_1\to \FF_2$ is a morphism in ${\bf Alg}_1$ and $\psi:\GG_1\to \GG_2$ is a group homomorphism such that $\psi\circ \iota_1=\iota_2\circ \varphi|_{\FF_1^\times}$.

\end{definition}

In particular, one has a natural (forgetful) projection functor $\pi:{\bf DecGr}\to {\bf Alg}_1$ by $\pi(\FF,G)=\FF$ and
$\pi(\varphi,\psi)=\varphi$.
Note also that for any object $(\FF,\GG,{\mathcal A})$ of ${\bf Gr}_\NN$ the pair $(\FF,\GG)$ is a decorated group, therefore, the projection $(\FF,\GG,{\mathcal A})\mapsto (\FF,\GG)$ defines a (forgetful) functor ${\bf Gr}_\NN\to {\bf DecGr}$.

\begin{definition} A {\it generalized $K_1$-theory} is a functor $K:\NN\to {\bf DecGr}$  such that $\pi\circ K=\pi|_\NN$.

\end{definition}

In what follows we will construct a number of generalized $K_1$-theories as  compositions of  an $\NN$-current group ${\mathfrak s}:\NN\to {\bf Gr}_\NN$ with a certain functors $\kappa$ from ${\bf Gr}_\NN$ to the category of {\it decorated groups}.

Given  an object $(\FF,\GG,{\mathcal A})$ of ${\bf Gr}_\NN$, we define $\kappa_{com}(\FF,\GG,{\mathcal A}):=\GG/[\GG,\GG]$, where
$[\GG,\GG]$ is the (normal) commutator subgroup of $G$.

\begin{lemma}
\label{le:commutative K1}
The correspondence $(\FF,\GG,{\mathcal A})\mapsto \kappa_{com}(\FF,\GG,{\mathcal A})$ defines a functor
$\kappa_{com}$ from ${\bf Gr}_\NN$ to the category of decorated abelian groups. In particular, for any
$\NN$-current group ${\mathfrak G}:\NN\to {\bf Gr}_\NN$ the composition $\kappa_{com}\circ {\mathfrak G}$ is a generalized $K_1$-theory.
\end{lemma}

Note that each $K_1$-theory defined by Lemma \ref{le:commutative K1} is still commutative. Below we will construct a number of non-commutative {\it nilpotent} $K_1$-theories in a similar manner.

\begin{definition}
\label{def:stable nilpotent} Let ${\mathcal A}$ be an
$\FF$-algebra. Recall that a subset $S$ of  ${\mathcal A}$ is
called a  nilpotent if  $S^n=0$ for some $n$. In particular, an
element $e\in {\mathcal A}$ is nilpotent if $\{e\}$ is nilpotent,
i.e., $e^n=0$. We denote by ${\mathcal A}_{nil}$ the set of all
nilpotent elements in ${\mathcal A}$.  Note that ${\mathcal
A}_{nil}$ is invariant under adjoint action of the group of units
${\mathcal A}^\times$. We say that an element $e\in {\mathcal A}$
is an {\it $\FF$-stable nilpotent} if $(\FF\cdot e\cdot \FF)^n=0$
for some $n> 0$. Denote by $(\FF,{\mathcal A})_{nil}$ the set of
all $\FF$-stable nilpotent element of ${\mathcal A}$. Denote also
by ${\mathcal A}_{nil}^{\FF^\times}$ and $(\FF,{\mathcal
A})_{nil}^{\FF^\times}$ the centralizers of $\FF^\times$ in
${\mathcal A}_{nil}$ and $(\FF,{\mathcal A})_{nil}$ respectively.
\end{definition}

In particular, taking ${\mathcal A}=\FF$, we see that $f\in (\FF,\FF)_{nil}$ if and only if the ideal $\FF f\FF\subset \FF$ is nilpotent.
Note also that if  the image $\iota(\FF)$ is in the center of ${\mathcal A}$, then  $(\FF,{\mathcal A})_{nil}={\mathcal A}_{nil}$.

Let $(\FF,\GG,{\mathcal A})$ be an object of ${\bf Gr}_\NN$ and  let $S$ be a subset of ${\mathcal A}_{nil}$, denote by  $E_S=E_S(\FF,\GG,{\mathcal A})$ the subgroup of $\GG$ generated by all $1+xsx^{-1}$, $s\in S$, $x\in \GG$. Clearly, $E_S$ is a normal subgroup of $\GG$. Then denote the quotient group $\GG/E_S$ by:
\begin{equation}
\label{eq:nilpotents and stable nilpotents}
\begin{cases} \kappa_{nil}(\FF,\GG,{\mathcal A}) & \text{if $S={\mathcal A}_{nil}$}\\
\kappa_{stnil}(\FF,\GG,{\mathcal A})& \text{if $S=(\FF,{\mathcal A})_{nil}$}\\
\kappa_{nil,inv}(\FF,\GG,{\mathcal A})& \text{if $S={\mathcal A}_{nil}^{\FF^\times}$}\\
\kappa_{stnil,inv}(\FF,\GG,{\mathcal A})& \text{if $S=(\FF,{\mathcal A})_{nil}^{\FF^\times}$}\\
\end{cases}
\end{equation}

\begin{lemma}
Each of the four correspondences
$$(\FF,{\mathcal L},{\mathcal A})\mapsto \kappa(\FF,{\mathcal L},{\mathcal A}) \ ,$$
where $\kappa=\kappa_{nil}$, $\kappa_{stnil}$, $\kappa_{nil,inv}$, $\kappa_{stnil,inv}$, defines a functor
$$\kappa:{\bf Gr}_\NN\to {\bf DecGr} \ .$$
In particular, for any $\NN$-current group ${\mathfrak G}:\NN\to {\bf Gr}_\NN$ the composition $\kappa\circ {\mathfrak G}$ is a
generalized
$K_1$-theory $\NN\to {\bf DecGr}$.
\end{lemma}

\begin{proof} Clearly, in each case of \eqref{eq:nilpotents and stable nilpotents}, the association ${\mathcal A}\mapsto S=S_{\mathcal A}$ is functorial, i.e., commutes with morphisms of  $\FF$-algebras. Therefore, the association ${\mathcal A}\mapsto E_S$ is also functorial in all four cases
\eqref{eq:nilpotents and stable nilpotents}. This finishes the proof of the lemma.
\end{proof}

We will elaborate examples of generalized (non-commutative) $K_1$-groups in a separate paper.

\subsection{$\NN$-current groups for compatible pairs}
\label{subsect:NN-current groups for compatible pairs}
Here we keep the notation of Section \ref{sect:Noncommutative current Lie algebras}. Since both  $\FF$ and $A$ are now  unital
algebras, so is  $\FF\otimes A=\FF\cdot A$.

Note that for the $\NN$-Lie algebra ${\mathfrak s}:\FF\mapsto (\FF,(\gg,A)(\FF),\FF\cdot A)$ the corresponding $\NN$-current group is of the form $(\FF,G_{\gg,A}(\FF),\FF\cdot A)$, where $G_{\gg,A}(\FF)$ is the normalizer of $(\gg,A)(\FF)$ in $(\FF\cdot A)^\times$ (i.e., $G_{\gg,A}(\FF)=\{g\in (\FF\cdot A)^\times:g\cdot (\gg,A)(\FF)\cdot g^{-1}=(\gg,A)(\FF)\}$).

%On the other hand, Lemma \ref{le:N-Lie correspondence} implies that the $\NN$-Lie algebra corresponding to the above $\NN$-group takes the form: $(\FF,(\gg,A)_N(\FF),\FF\times A)$, where  %$(\gg,A)_N(\FF)$ is the normalizer  Lie algebra of $(\gg,A)(\FF)$ in $\FF\cdot A$.

The following facts are obvious.

\begin{lemma}
\label{le:general stabilization} Let $(\gg,A)$ be a compatible
pair and $S\subset \gg$ be a generating set of $\gg$ as a Lie algebra. Then
for any object $\FF$ of ${\bf Alg}_1$ an element $g\in (\FF\cdot
A)^\times$ belongs to $G_{\gg,A}(\FF)$ if and only if:
\begin{equation}
\label{eq:general stabilization}
g(u\cdot  x) g^{-1}\subset (\gg,A)(\FF)
\end{equation}
for all $x\in S$, $u\in \FF$.

\end{lemma}

\begin{lemma} For each compatible pair of the form  $(\gg,A)=(sl_n(\kk),M_n(\kk))$ and an object $\FF$ of $\NN_1$ one has: $G_{\gg,A}(\FF)=GL_n(\FF)=(\FF\cdot A)^\times$.

%\noindent (a) $g\cdot (\gg,A)_N(\FF)\cdot g^{-1}=(\gg,A)_N(\FF)$ for all $g\in G_{\gg,A}(\FF)$.

%and $(\gg,A)_N(\FF)=gl_n(\FF)=\FF\cdot A$.

\end{lemma}

In what follows we will consider compatible pairs of the form $(\gg,End(V))$ where $V$ is a simple finite-dimensional $\gg$-module. By choosing an appropriate basis in $V$, we identify $A=End(V)$ with $M_n(\kk)$ so that $G_{\gg,A}(\FF)\subset GL_n(\FF)$. In all cases to be considered, we will compute the ``Cartan subgroup'' $(\FF^\times)^n\cap G_{\gg,A}(\FF)$ of $G_{\gg,A}(\FF)$.

Let  $\Phi_0$ be the  bilinear form on the $\kk$-vector space
$V=\kk^n$ given by:
$$\Phi_0(x, y)=x_1y_n+x_2y_{n-1}+\dots +x_ny_1.$$

Also define the bilinear form $\Phi _1$ on $\kk^{2m}$ by:
$$\Phi_1(x, y)=x_1y_{2m}+x_2y_{2m-1}+\dots +x_my_m-x_{m+1}y_{m-1}\dots -x_{2m}y_1.$$

\begin{proposition} 
\label{pr:cartan for classical groups} Let $\FF$ be an object of ${\bf Alg}_1$,
$A=M_n(\kk)$, and suppose that either $\gg=o(\Phi _0)$ or $\gg=o(\Phi _1)$ and $n=2m$.
Then an invertible diagonal matrix
$D=diag(f_1,...,f_n)\in GL_n(\FF)$ belongs to $G_{\gg,A}(\FF)$ if and only if
$$f_if_{n-i+1}-f_1f_n\in I_1(\FF)=\FF[\FF,\FF]$$
for $i=1,\ldots, n$.
\end{proposition}
\begin{proof} We will prove the proposition for  $\gg=o(\Phi_0)$ (the proof for $\gg=o(\Phi _1)$ is nearly identical). It is easy to see
that  $\gg$ is  a Lie subalgebra of $sl_n(\kk)$  generated by all $e_{ij}:=E_{ij}-E_{n-j+1,
n-i+1}$, $i,j=1,\ldots,n$. therefore, Lemma \ref{le:general stabilization} (with $S=\{e_{ij}\}$) guarantees
that $D=(f_1,\ldots,f_n)\in (\FF^\times)^n$
belongs to $G_{\gg,A}(\FF)$ if and only if $D(u\cdot e_{ij})D^{-1}\subset (\gg, A)(\FF)$ for all $u\in \FF$, $i,j=1,\ldots,n$.
Note that
$$D(u\cdot e_{ij})D^{-1}=f_iuf_j^{-1}E_{ij}-f_{j'}uf_{i'}^{-1}E_{j', i'}=f_iuf_j^{-1}e_{ij}+\delta_{ij}(u)E_{j', i'} \ ,$$
where $\delta_{ij}(u)=f_iuf_j^{-1}-f_{j'}uf_{i'}^{-1}$ and
$i'=n+1-i$, $j'=n+1-j$. Therefore,  taking into account that
$(\gg, A)(\FF)=[\FF,\FF]\cdot 1+\FF \cdot \gg+I_1(\FF)\cdot
sl_n(\kk)$ by Corollary \ref{cor:orthogonal and simplectic lie
algebra}, we see that $D(u\cdot e_{ij})D^{-1}\in (\gg, A)(\FF)$ if and only if $\delta_{ij}(u)\in I_1(\FF)$. Note
that $\delta_{ij}(u)\equiv u\delta_{ij}(1) \mod I_1(\FF)$. Since
$I_1(\FF)$ is an ideal, $u\delta_{ij}(1)\in I_1(\FF)$ for all
$u\in \FF$ if and only if $\delta_{ij}(1)\in I_1(\FF)$, i.e.,
$f_if_j^{-1}- f_{j'}f_{i'}^{-1}\mod I_1(\FF)$. Taking into
account that $f_if_j^{-1}- f_{j'}f_{i'}^{-1}\equiv
(f_if_{i'}-f_jf_{j'})f_j^{-1}f_{i'}^{-1}\mod I_1(\FF)$, we see
that $D\in G_{\gg,A}(\FF)$ if and only if
$f_if_{n+1-i}-f_jf_{n+1-j}\in I_1(\FF)$ for all $i, j$. Clearly,
it suffices to take $j=1$. The proposition is proved.
\end {proof}

%Define now two more subgroups $E(\gg,A,\FF)$ and $\overline{E(\gg,A,\FF)}$
%of $(\FF\cdot A)^\times$ as follows.
%Clearly, for each nilpotent element $x\in \FF\cdot A$  the element $1+x$ belongs to $(\FF\cdot A)^\times$.

%\begin{proposition} For each compatible pair $(\gg,A)$ and an object $\FF$ of $\NN_1$ one has:
%
%\noindent (a)  $E(\gg,A,\FF)$ is a subgroup of $E_N(\gg,A,\FF)$.
%
%\noindent (b) Both $E(\gg,A,\FF)$ and $E_N(\gg,A,\FF)$ are normal subgroups of  $G(\gg,A,\FF)$.
%\end{proposition}
%%%%%%%%%%%%%%%%%%%%%%%%%
To formulate the main result of this section we need the following notation.

For any $\ell\ge 0$ and any $m_1,\ldots,m_{\ell+1}\in \FF$ denote
\begin{equation}
\label{eq:k-th divided difference}
\Delta^{(\ell)}(m_1,\ldots,m_{\ell+1})=\sum_{k=0}^\ell (-1)^k \binom{\ell}{k} m_{k+1}
\end{equation}
and refer to it as the $\ell$-th difference derivative. Clearly,
$$\Delta^{(\ell)}(m_1,\ldots,m_{\ell+1})=\Delta^{(\ell-1)}(m_1,\ldots,m_\ell)-\Delta^{(\ell-1)}(m_2,\ldots,m_{\ell+1}) \ .$$

Let $A_n=M_n(\kk)=End(V_{n-1})$, where
$V_{n-1}$ is the $n$-dimensional irreducible  $sl_2(\kk)$-module. Then $G_{sl_2(\kk),A_n}(\FF)$ is naturally a subgroup of $GL_n(\FF)$.

\begin{maintheorem}
\label{th:sl2 in sln}
For any  object $\FF$ of  ${\bf Alg}_1$,  the ``Cartan subgroup'' $(\FF^\times)^n \cap G_{sl_2(\kk),A_n}(\FF)$
consists of all $D=(f_1,...,f_n)\in (\FF^\times)^n$ such that:
\begin{equation}
\label{eq:difference derivative}
\Delta^{(k)}(f_1f_2^{-1},\ldots,f_{k+1} f_{k+2}^{-1}) \in I_k(\FF)
\end{equation}
%\begin{equation}
%\label{eq:difference derivative2}
%\Delta^{(k)}(m_{n-1}^{-1},\ldots,m_{n-k-1}^{-1})\in I_k(\FF)
%\end{equation}
for $k=1,\ldots,n-2$.
%, where  $m_i=f_i f_{i+1}^{-1}$, $i=1,...,n-1$.
\end{maintheorem}

\begin{proof} We will prove the theorem in several steps. 
First, we prove Proposition \ref{pr:reformulation of D-conjugation}
by using Lemmas \ref{le:above diagonal sl2 in sln} and \ref{le:D-conjugation of E}. Then we prove Lemma
\ref{le:mij} and Proposition \ref{pr:mij and mij*}. The proof of the Proposition is based on Lemma \ref{le:tr-in}.
The final step in the proof of Theorem \ref{th:sl2 in sln} is 
Theorem \ref{th:difference derivative implications}. This Theorem required Proposition \ref{pr:dij homogeneous}, Lemma
\ref{le: Delta(j-i)}, and Proposition \ref{pr:degree k noncommutative operator}. Our proofs of Propositons
\ref{pr:dij homogeneous} and \ref{pr:degree k noncommutative operator} use Lemmas \ref{le:expl dij} and
\ref{le:partial decomposition} correspondingly. 

We start with a
characterization of $(\FF^\times)^n\cap G_{sl_2(\kk),A_n}(\FF)$.
\begin{proposition}
\label{pr:reformulation of D-conjugation}
A diagonal matrix $D=(f_1,\ldots,f_n)\in (\FF^\times)^n$ belongs to the group $G_{sl_2(\kk),A_n}(\FF)$ if and only if:
\begin{equation}
\label{eq:particular stabilization1}
\Delta^{(k)}(f_1uf_2^{-1},\ldots,f_{k+1}uf_{k+2}^{-1})\in
I_k(\FF),
\end{equation}
\begin{equation}
\label{eq:particular stabilization2}
\Delta^{(k)}(f_n uf_{n-1}^{-1},\ldots,f_{n-k} uf_{n-1-k}^{-1})\in I_k(\FF) \ .
\end{equation}
for  $k=1,\ldots, n-2$ and all $u\in \FF$.

\end{proposition}

\begin{proof} 
%Denote by $G_0$ the set of all $D\in (\FF^\times)^{\ell+1}\subset GL_{\ell+1}(\FF)$ such that 
%$D\cdot E \cdot D^{-1}\in (sl_2(\kk),A_\ell)(\FF)$ and $D\cdot F \cdot D^{-1}\in (sl_2(\kk),A_\ell)(\FF)$. 
%All we have to prove is that $G_0\subset  G(sl_2(\kk),A_\ell,\FF)$.
%Note that $A_n=V_0+V_2+\cdots + V_{2(n-1)}$in the notation of Theorem \ref{th:perfect sl_2} and the sum is direct.
Denote by $(A_n)_k$ the set of all $x\in A_n$ such that
$[H,x]=kx$. Clearly, $(A_n)_k\ne 0$ if and only if $k$ is even and
$-2(n-1)\le k\le 2(n-1)$. In fact, $(A_n)_{2k}$ is the span of
all those $E_{ij}$ such that $2(j-i)=k$. In particular, $E\in
(A_n)_2$. Denote also
$(sl_2(\kk),A_n)(\FF)_k:=(sl_2(\kk),A_n)(\FF)\cap \FF\cdot
(A_n)_k$.

\begin{lemma}
\label{le:above diagonal sl2 in sln} 
The components $(sl_2(\kk),A_n)(\FF)_j$, $j=-2, 2$ are given by:
\begin{equation}
\label{eq:above diagonal sl2 in sln}
(sl_2(\kk),A_n)(\FF)_2=\sum_{k=0}^{n-2} I_k(\FF)\cdot E^{(k)}, ~(sl_2(\kk),A_n)(\FF)_{-2}=\sum_{k=0}^{n-2} I_k(\FF)\cdot F^{(k)} \ ,
\end{equation}
where
$$E^{(k)}=\sum_{i=k+1}^{n-1} i\binom{i-1}{k}E_{i,i+1}, F^{(k)}=\sum_{i=k+1}^{n-1} i\binom{i-1}{k}E_{n+1-i,n-i}$$
for $k=0,1,\ldots, n-2$ form a basis for $(A_n)_2$ and $(A_n)_{-2}$ respectively.
\end{lemma}

\begin{proof} Let us prove the formula for $j=2$. Let $p_0(x), \ldots p_{n-2}(x)$ be any polynomials 
in $\kk[x]$ such that $\deg p_k(x)=k$ for all $k$. 
Then it is easy to see that
$$(sl_2(\kk),A_n)(\FF)_2=\sum_{k=0}^{n-2} I_k(\FF)\cdot ( p_k(H) \cdot E) \ .$$
Take $p_k(H)=\binom{H'}{k}$, where $H'=\frac{1}{2}(n\cdot 1-H)=\sum_{i=1}^n  (i-1) E_{ii}$. 
Then $p_k(H)=\sum_{i=1}^n \binom{i-1}{k}E_{ii}$ and
$$p_k(H) \cdot E=\left(\sum_{i=1}^n \binom{i-1}{k}E_{ii}\right)\left(\sum_{i=1}^n iE_{i,i+1}\right)=E^{(k)} \ .$$
To prove the formula for $j=-2$, it suffices to conjugate the formula for $j=2$ with 
the matrix of the longest permutation $w_0=\sum_{i=1}^n E_{i,n+1-i}$, i.e., apply the involution $E_{ij}\mapsto E_{n+1-i,n+1-j}$.
\end{proof}

\begin{lemma}
\label{le:D-conjugation of E}
For each diagonal matrix $D=(f_1,f_2,\ldots,f_n)\in (\FF^\times)^n$ and $u\in \FF$ one has
$$D(u\cdot E)D^{-1}=\sum_{k=0}^{n-2} \Delta^{(k)}(f_1uf_2^{-1},\ldots,f_{k+1}uf_{k+2}^{-1})\cdot (-1)^kE^{(k)} \ ,$$
$$D(u\cdot F)D^{-1}=\sum_{k=0}^{n-2} \Delta^{(k)}(f_n uf_{n-1}^{-1},\ldots,f_{n-k} uf_{n-1-k}^{-1})\cdot (-1)^{k+1}F^{(k)} \ ,$$
where  $\Delta^{(k)}$ is the $k$-th divided difference
as in  \eqref{eq:k-th divided difference}.
\end{lemma}

\begin{proof} It is easy to see that the elements $E^{(k)}, F^{(k)}$ satisfy:
\begin{equation}
\label{eq:Eii+1}
iE_{i,i+1}=\sum_{k=i-1}^{n-2} (-1)^{k+1-i} \binom{k}{i-1} E^{(k)}, 
~iE_{n+1-i,n-i}=\sum_{k=i-1}^{n-2} (-1)^{k+1-i} \binom{k}{i-1} F^{(k)}
\end{equation}
for $i=1,\ldots,n-1$.

Furthermore,
$$D(u\cdot E)D^{-1}=\sum_{i=1}^{n-1} f_iuf_{i+1}^{-1}\cdot iE_{i,i+1}=
\sum_{i=1}^{n-1} f_iuf_{i+1}^{-1}\cdot \sum_{k=0}^{n-2} (-1)^{k+1-i} \binom{k}{i-1} E^{(k)}=$$
$$\sum_{k=0}^{n-2} \sum_{i=1}^{n-1} (-1)^i \binom{k}{i-1}f_iuf_{i+1}^{-1}\cdot (-1)^{k+1}E^{(k)}$$
$$=\sum_{k=i-1}^{n-2}\Delta
^{(k)}(f_1uf_2^{-1},\ldots,f_{k+1}uf_{k+2}^{-1})\cdot
(-1)^{k+1}E^{(k)} \ .$$ The formula for $D(u\cdot F)D^{-1}$
follows. The lemma is proved.
\end{proof}

Now we are ready to finish the proof of Proposition \ref{pr:reformulation of
D-conjugation}.

Since the set $S=\{E,F\}$ generates $sl_2(\kk)$,
Lemma \ref{le:general stabilization} guarantees that $D\in
GL_n(\FF)$ belongs to $G_{sl_2(\kk),A_n}(\FF)$ if and only if
$D(u\cdot E)D^{-1}, D(u\cdot F)D^{-1}\in (sl_2(k),A_n)(\FF)$ for
all $u\in \FF$. Using  the obvious fact that $D(u\cdot
E)D^{-1}\subset \FF\cdot (A_n)_2$ for all $D\in (\FF^\times)^n$,
$u\in \FF$,  we see that $D(u\cdot E)D^{-1}\in
(sl_2(\kk),A_n)(\FF)$ if and only if $D(u\cdot E)D^{-1}\in
(sl_2(\kk),A_n)(\FF)_2$. In turn, using Lemmas \ref{le:above
diagonal sl2 in sln} and \ref{le:D-conjugation of E}, we see that
this is equivalent to
$$D(u\cdot E)D^{-1}=\sum_{k=0}^{n-2} \Delta^{(k)}(m_1,\ldots,m_{k+1})\cdot (-1)^{k+1}E^{(k)}\in \sum_{k=0}^{n-2} I_k(\FF)\cdot E^{(k)} \ ,$$
which, because the $E^{(0)},\ldots,E^{(n-2)}$ are linearly
independent, is equivalent to \eqref{eq:particular
stabilization1}. Applying the above argument to $D(u\cdot
F)D^{-1}$, we obtain
$$D(u\cdot E)D^{-1}=\sum_{k=0}^{n-2} \Delta^{(k)}(f_n uf_{n-1}^{-1},\ldots,f_{n-k} uf_{n-1-k}^{-1})
\cdot (-1)^{k+1}F^{(k)}\in \sum_{k=0}^{n-2} I_k(\FF)\cdot F^{(k)} \ ,$$
which gives \eqref{eq:particular stabilization2}.

The proposition is proved.
\end{proof}

Furthermore, we  need to establish
some basic properties of inclusions \eqref{eq:difference
derivative}.

\begin{lemma}
\label{le:mij} 
Let $m_1,m_2,\ldots,m_\ell$ be elements of $\FF$. The following are equivalent:

\noindent (a) $\Delta^{(k)}(m_1,\ldots,m_{k+1})\in I_k(\FF)$ for
all $1\le k\le \ell-1$.

\noindent (b) $\Delta^{(j-i)}(m_i,\ldots,m_j)\in I_{j-i}(\FF)$ for all $1\le i\le j\le \ell$.

\end{lemma}

\begin{proof} The implication (b) $\Rightarrow $(a) is obvious. Prove the implication (a)$\Rightarrow $(b). Denote
\begin{equation}
\label{eq:mij}
m_{ij}:=\Delta^{(j-i)}(m_i,\ldots,m_j)
\end{equation}
for all $1\le i\le j\le \ell$.  In particular, $m_{ii}=m_i$ and:
\begin{equation}
\label{eq:recursive difference derivative}
m_{ij}=m_{i,j-1}-m_{i+1,j}
\end{equation}
for all $1\le i\le j\le \ell$.

Prove that the inclusions $m_{1,k+1}\in I_k(\FF)$ for $0\le k\le
\ell$ imply inclusions $m_{ij}\in I_{j-i}(\FF)$ for all $i\le j$
such that $j-i\le \ell$. We proceed by induction on $i$. The
basis of the induction, when $i=1$, is obvious. Assume that $i>1$
and $j\le \ell-1-i$. Then the inclusions (the inductive
hypothesis)
$$m_{i-1,j}=m_{i-1,j-1}-m_{ij}\in I_{j+1-i}(\FF)\subset I_{j-i}(\FF)$$
and $m_{i-1,j-1}\in I_{j-i}(\FF)$ imply that $m_{ij}\in
I_{j-i}(\FF)$. This finishes the proof of the implication
(a)$\Rightarrow $(b). The lemma is proved.
\end{proof}

\begin{proposition}
\label{pr:mij and mij*} 
Let $m_1,m_2,\ldots,m_\ell$ be invertible
elements of $\FF$. The following  are equivalent:

\noindent (a) $\Delta^{(j-i)}(m_i,\ldots,m_j)\in I_{j-i}(\FF)$ for all $i\le j$.

\noindent (b) $\Delta^{(j-i)}(m_i^{-1},\ldots,m_j^{-1})\in I_{j-i}(\FF)$ for all $i\le j$.

\end{proposition}

\begin{proof} We need the following notation.
Similarly to \eqref{eq:mij} denote
\begin{equation}
\label{eq:mij and mij*}
m_{ij}:=\Delta^{(j-i)}(m_i,\ldots,m_j),~m_{ij}^*:=\Delta^{(j-i)}(m_i^{-1},\ldots,m_j^{-1})
\end{equation}
for all $1\le i\le j\le \ell$. In particular, $m_{ii}^*=m_i^{-1}$
and $m_{12}^*=-m_1^{-1}(m_1-m_2)m_2^{-1}=-m_{11}^*m_{12}m_{22}^*$.

We need  the following recursive formula for $m_{ij}^*$.

\begin{lemma} 
\label{le:tr-in} 
In the above notation, we have for all $1\le i<j\le \ell$:
\begin{equation}
\label{eq:mij* recursive}
m_{ij}^*=\doublesubscript{\sum}{i\le i_1\le j_1\le j,i\le i_2<j_2\le j,i\le i_3\le j_3\le j}{j_1-i_1+j_2-i_2+j_3-i_3=j-i} 
c_{i,i_1,i_2,i_3}^{j,j_1,j_2,j_3} m^*_{i_1,j_1}m_{i_2,j_2} m^*_{i_3,j_3} \ ,
\end{equation}
where the coefficients are translation-invariant integers:
$$c_{i+1,i_1+1,i_2+1,i_3+1}^{j+1,j_1+1,j_2+1,j_3+1}=c_{i,i_1,i_2,i_3}^{j,j_1,j_2,j_3} \ .$$
%m_{kk}^*m_{ij}*m_{\ell,\ell}^*+
\end{lemma}

\begin{proof} We proceed by induction on $j-i$. If $j=i+1$, we obtain:
$$m_{i,i+1}^*=-m_i^{-1}(m_i-m_{i+1})m_{i+1}^{-1}=-m_{ii}^*m_{i,i+1}m_{i+1,i+1}^*=-m_{i+1,i+1}^*m_{i,i+1}m_{i,i}^* \ .$$
Next, assume that $j-i>1$. Then, using the translation invariance
of the coefficients in \eqref{eq:mij* recursive} for
$m_{i,j-1}^*$, we obtain
$$m_{ij}^*=m_{i,j-1}^*-m_{i+1,j}^*=\doublesubscript{\sum}{i\le i_1\le j_1\le j-1,i\le i_2<j_2\le j-1,i\le i_3
\le j_3\le j-1}{j_1-i_1+j_2-i_2+j_3-i_3=j-1-i} c_{i,i_1,i_2,i_3}^{j-1,j_1,j_2,j_3} \delta_{i_1,i_2,i_3}^{j_1,j_2,j_3} \ , $$
where $\delta_{i_1,i_2,i_3}^{j_1,j_2,j_3}=m^*_{i_1,j_1}m_{i_2,j_2} m^*_{i_3,j_3}-m^*_{i_1+1,j_1+1}m_{i_2+1,j_2+1} m^*_{i_3+1,j_3+1}$.
Furthermore,
$$\delta_{i_1,i_2,i_3}^{j_1,j_2,j_3}=m^*_{i_1,j_1+1}m_{i_2,j_2} m^*_{i_3,j_3}+
m^*_{i_1+1,j_1+1}(m_{i_2,j_2} m^*_{i_3,j_3}-m_{i_2+1,j_2+1}m^*_{i_3+1,j_3+1})$$
$$=m^*_{i_1,j_1+1}m_{i_2,j_2} m^*_{i_3,j_3}+m^*_{i_1+1,j_1+1}m_{i_2,j_2+1} m^*_{i_3,j_3} +
m^*_{i_1+1,j_1+1}m_{i_2+1,j_2+1}m^*_{i_3,j_3+1} \ .$$
This proves the formula \eqref{eq:mij* recursive} for $m_{ij}^*$. The lemma is proved.
\end{proof}

\noindent We are ready to finish the proof of  Proposition  \ref{pr:mij and mij*} now.

Due to the symmetry, it suffices to prove only one implication, say (a)$\Rightarrow $(b). The desired
implication follows inductively from \eqref{eq:mij and mij*}.
\end{proof}

Note that, in the view of Lemma \ref{le:mij}, the condition (a)
(resp. the condition (b)) of Proposition \ref{pr:mij and mij*} for
$m_i=f_if_{i+1}^{-1}$, $i=1,\ldots,n-2$, is a particular case of
\eqref{eq:particular stabilization1} (resp. of
\eqref{eq:particular stabilization2}) with $u=1$. 
%Therefore,  to finish the proof of Theorem \ref{th:sl2 in sln}, we need to prove
%the converse. Since the conditions (a) and (b) of Proposition
%\ref{pr:mij and mij*} are equivalent, all we have to do is to
%prove the following result.

Furthermore, we need one more result in order to finish the proof of Theorem \ref{th:sl2 in sln}.

\begin{theorem}
\label{th:difference derivative implications}
The inclusions \eqref{eq:difference derivative} imply the inclusions \eqref{eq:particular stabilization1} (and, therefore,
\eqref{eq:particular stabilization2}).

\end{theorem}

\begin{proof}  To prove the theorem we will develop a formalism of homogeneous
maps $\FF\to \FF$ (relative to the ideals $I_k(\FF))$.

\begin{definition}
We say that a $\kk$-linear map $\partial:\FF\to \FF$ is {\it homogeneous of degree} 
$\ell$ if $\partial(I_k(\FF))\subset I_{k+\ell}(\FF)$ for all $k\ge 0$;  denote by $End^{(\ell)}(\FF)$ the set of all such maps.

\end{definition}

Lemma \ref{le:ideals I_k} guarantees that  for each $f_1\in I_{\ell_1}(\FF)$, 
$f_2\in I_{\ell_2}(\FF)$, the map $\FF\to \FF$ given by $u\mapsto f_1uf_2$ is homogeneous of degree $\ell_1+\ell_2$.

We construct a number of homogeneous maps of degree $1$ as follows. For an invertible element $m\in \FF^\times$ define $\partial_m:\FF\to \FF$ by
$$\partial_m(u)=m u m^{-1}-u =[m,u m^{-1}]\ .$$
Clearly, $\partial_m:\FF\to \FF$ is homogeneous of degree $1$.

\begin{proposition}
\label{pr:dij homogeneous} 
Let $m_1,m_2,\ldots,m_\ell$ be invertible elements of $\FF$ such that, in the notation \eqref{eq:mij}, one has
$m_{ij}\in I_{j-i}(\FF)$  for all $1\le i\le j \le \ell$.
 Then
$$\Delta^{(j-i)}(\partial_{m_i},\ldots,\partial_{m_j})\in End^{(j+1-i)}(\FF)$$
for all $1\le i\le j\le \ell$.
\end{proposition}

\begin{proof} We need the following notation.
Similarly to \eqref{eq:mij}, denote
\begin{equation}
\label{eq:dij}
\partial_{ij}=\Delta^{(j-i)}(\partial_{m_i},\ldots,\partial_{m_j})
\end{equation}
for $1\le i\le j\le \ell$. By definition, $\partial_{ii}=\partial_{m_i}$ and $\partial_{i,i+1}=\partial_{m_i}-\partial_{m_{i+1}}$.

\begin{lemma} 
\label{le:expl dij}
For each $u\in \FF$ and $1\le i\le j\le \ell$ one has:
\begin{equation}
\label{eq:explicit dij}
\partial_{ij}(u)=\doublesubscript{\sum}{i\le i_1\le j_1\le j,i\le i_2\le j_2\le j}
{j_1-i_1+j_2-i_2=j-i} c_{i,i_1,i_2}^{j,j_1,j_2} [m_{i_1,j_1},u m^*_{i_2,j_2}]
\end{equation}
in the notation \eqref{eq:mij and mij*}, where the coefficients are translation-invariant integers:
$$c_{i+1,i_1+1,i_2+1}^{j+1,j_1+1,j_2+1}=c_{i,i_1,i_2}^{j,j_1,j_2} \ .$$
%m_{kk}^*m_{ij}*m_{\ell,\ell}^*+
\end{lemma}

\begin{proof} We proceed by induction on $j-i$. If $j=i$, we have $\partial_{ii}(u)=\partial_{m_i}(u)=[m_i,um_i^{-1}]$.
%$$m_{i,i+1}= -m_i^{-1}(m_i-m_{i+1})m_{i+1}^{-1}=-m_{ii}^*m_{i,i+1}m_{i+1,i+1}^*=-m_{i+1,i+1}^*m_{i,i+1}m_{i,i}^* \ .$$
Next, assume that $j-i>0$. Then, using the translation-invariance
of the coefficients in \eqref{eq:explicit dij} for
$\partial_{i,j-1}(u)$, we obtain
$$\partial_{ij}(u)=\partial_{i,j-1}(u)-\partial_{i+1,j}(u)=\doublesubscript{\sum}{i\le i_1\le j_1\le j-1,i\le i_2<j_2\le j-1}{j_1-i_1+j_2-i_2=j-1-i} c_{i,i_1,i_2}^{j-1,j_1,j_2} \delta_{i_1,i_2}^{j_1,j_2} \ , $$
where $\delta_{i_1,i_2}^{j_1,j_2}=[m_{i_1,j_1},u m^*_{i_2,j_2}]-[m_{i_1+1,j_1+1},u m^*_{i_2+1,j_2+1}]$.
Furthermore,
$$\delta_{i_1,i_2}^{j_1,j_2}=[m_{i_1,j_1+1},u m^*_{i_2,j_2}]+[m_{i_1+1,j_1+1},u m^*_{i_2,j_2}]-[m_{i_1+1,j_1+1},u m^*_{i_2+1,j_2+1}]$$
$$=[m_{i_1,j_1+1},u m^*_{i_2,j_2}]+[m_{i_1+1,j_1+1},u m^*_{i_2,j_2+1}]\ .$$
This proves the formula \eqref{eq:explicit dij} for $\partial_{ij}(u)$. The lemma is proved.
\end{proof}

Now we can finish the proof of Proposition \ref{pr:dij homogeneous}.

 Since $m_{i_1,j_1}\in I_{j_1-i_1}(\FF)$, $m_{i_2,j_2}^*\in
I_{j_2-i_2}(\FF)$, and $[m_{i_1,j_1},u m^*_{i_2,j_2}]$ belongs to
the ideal $I_{k+j_1-i_1+j_2-i_2+1}(\FF)$ for all $u\in I_k(\FF)$,
formula \eqref{eq:explicit dij} guarantees inclusion
$\partial_{ij}(I_k(\FF))\subset I_{k+j+1-i}(\FF)$ for all $k\ge
0$. This proves  Proposition \ref{pr:dij homogeneous}.
\end{proof}

We continue to study ``difference operators" in $\FF$. Let
$m_i,d_i:\FF\to \FF$, $i=1,\ldots,\ell$ be linear maps. Denote
\begin{equation}
\label{le:upper partial}
\partial^{(i,j)}:=\Delta^{(j-i)}(m_id_{i+1}\cdots d_j,m_{i+1}d_{i+2}\cdots d_j,\ldots,m_{j-1}d_j,m_j)
\end{equation}
%where $m_{ij}:=\Delta^{(i-j)}(m_i,\ldots,m_j),$\partial_{ij}=\Delta^{(i-j)}(\partial_i,\ldots,\partial_j)$
for all $1\le i\le j\le \ell$.
%Given invertible elements $m_1,m_2,\ldots,m_\ell$  of $\FF$, define the maps $\partial^{(i,j)}:\FF\to \FF$, $1\le i\le j\le \ell$ by the formula:
%$$\partial^{(i,j)}=\sum_{k=0}^{j-i} (-1)^k \binom{j-i}{k} m_{k+i} \prod_{p=k+i}^j (\partial_{m_p}+1)$$
%for $1\le i\le j\le \ell$.

For instance, $\partial^{(i,i)}=m_i$, $\partial^{(i,i+1)}=m_id_{i+1}-m_{i+1}$, and
$$\partial^{(i,i+2)}=m_i d_{i+1}d_{i+2}-2m_{i+1}d_{i+2}+m_{i+2} \ .$$

\begin{lemma}
\label{le: Delta(j-i)} 
Let $D=(f_1,\ldots,f_n)\in (\FF^\times)^n$.
Set $m_i=f_if_{i+1}^{-1}$.  Then for each $u\in \FF$ one has
\begin{equation}
\label{eq:from delta to d}
\Delta^{(j-i)}(f_iuf_{i+1}^{-1},\ldots,f_juf_{j+1}^{-1})=\partial^{(i,j)}(u') \ ,
\end{equation}
where we abbreviated $u'=f_{j}uf_{j}^{-1}$ and $d_i:=\partial_{m_i}+1$.
\end{lemma}

\begin{proof} Indeed, for $i\le k\le j$ one has
$$f_kuf_{k+1}^{-1}=m_km_{k+1}\cdots m_j f_{j+1} u f_{j+1}^{-1}(m_{k+1}\cdots m_j)^{-1}=m_k(\partial_{m_{k+1}\cdots m_j}+1)(u')$$
$$=m_k(\partial_{m_k}+1)\cdots (\partial_{m_j}+1)(u') =m_kd_{k+1}\cdots d_j(u')\ .$$
Substituting so computed $f_kuf_{k+1}^{-1}$ into \eqref{eq:k-th divided difference}, we obtain \eqref{eq:from delta to d}. \end{proof}

Therefore, all we need to finish the proof of Theorem \ref{th:difference derivative implications} is to prove the following result.
\begin{proposition}
\label{pr:degree k noncommutative operator}
Let $m_1,m_2,\ldots,m_\ell$ be elements of $\FF$ and $\ell\ge 0$ such that, 
in the notation \eqref{eq:mij}, one has  $m_{ij}\in I_{j-i}(\FF)$ for all $1\le i\le j\le \ell$. Then  $$\partial^{(i,j)}\in End^{(j-i)}(\FF)$$
for all $1\le i\le j\le \ell$, where again we abbreviated $\partial_i:=\partial_{m_i}$.
\end{proposition}

\begin {proof} We need the following notations. Let
$M_I$ be the linear map $\FF\to \FF$ given by:
$$M_I=m_{i_0,i_1} \partial_{i_1+1,i_2} \partial_{i_2+1,i_3}\cdots \partial_{i_{k-1}+1,i_k} \ ,$$
where
$$m_{i',j'}:=\Delta^{(j'-i')}(m_{i'},\ldots,m_{j'}),~\partial_{i',j'}=
\begin{cases} d_{i'}-1 & \text{if $i'=j'$}\\
\Delta^{(j'-i')}(d_{i'},\ldots,d_{j'}) & \text{if $i'<j'$}
\end{cases}$$
for all $1\le i'\le j'\le \ell$.

\begin{lemma}
\label{le:partial decomposition} 
In the notation \eqref{le:upper partial} one has
\begin{equation}
\label{eq:upper partial I}
\partial^{(i,j)}=\sum_{I} c_I M_I \ ,
\end{equation}
where the summation is taken over all subsets
$I=\{i_0<i_1<i_2\cdots < i_k\}$ of $\{1,\ldots,\ell\}$ such that
$i_0=i$, $i_k=j$, and the coefficients $c_I\in \ZZ$ are
translation-invariant:
$$c_{I+1}= c_I\ , $$
where for any subset $I=\{i_0,\ldots,i_k\}\subset
\{1,\ldots,\ell\}$, we  abbreviate
$I+1=\{i_0+1,\ldots,i_k+1\}\subset \{2,\ldots,\ell+1\}$.
\end{lemma}

\begin{proof} We proceed by induction on $j-i$. The  basis of the induction when $j=i$ is obvious because
$\partial^{(i,i)}=M_I=m_i$ for all $i$, where $I=\{i\}$. Note that
$$\partial^{(i,j)}=\partial^{(i,j-1)}d_j-\partial^{(i+1,j)}=\partial^{(i,j-1)}\partial_{j,j}+(\partial^{(i,j-1)}-\partial^{(i+1,j)})$$
for all $i,j$. Therefore, we have by the inductive hypothesis and the translation-invariance of the coefficients $c_I$:
$$\partial^{(i,j)}=\sum_I c_I M_{I\sqcup\{j\}}+\sum_I c_I (M_I-M_{I+1}) \ ,$$
where the summations are over all subsets $I=\{i_0<i_1<\cdots< i_k\}$ of $\{1,\ldots,\ell\}$ such that
$i_0=i$, $i_k=j-1$ (we have used the fact that $M_I \partial_{j,j}=M_{I\sqcup\{j\}}$). It is easy to see that for any
$I=\{i_0<i_1<\cdots< i_k\}$ one has
$$M_I-M_{I+1}=M_{I_1}+M_{I_2}\cdots +M_{I_k} \ ,$$
where $I_j=\{i_0<i_1<i_2\cdots <i_{j-1}<i_j+1< \cdots <i_k+1\}$
for $j=1,2,\ldots,k$. Therefore,
$$\partial^{(i,j)}=\sum_I c_I M_{I\sqcup\{j\}}+\sum_{I,j} c_I M_{I_j} \ ,$$
i.e., $\partial^{(i,j)}$ is of the form \eqref{eq:upper partial I}. The lemma is proved.
\end{proof}

\noindent Now we can finish the proof Proposition \ref{pr:degree k noncommutative
operator}.

Recall from \eqref{eq:difference derivative} that
$m_{ij}=\Delta^{(j-i)}(m_i,\ldots,m_j)\in I_{j-i}(\FF)$ (hence
the map $u\mapsto m_{ij}u$ belongs to $End^{(j-i)}(\FF)$) and
from Proposition \ref{pr:dij homogeneous} that
$\partial_{ij}=\Delta^{(j-i)}(\partial_{m_i},\ldots,\partial_{m_j})\in
End^{(j+1-i)}(\FF)$ for all $1\le i\le j\le \ell$. This implies
that for $I=\{i_0<i_1<i_2<\cdots< i_k\}$ one has:
$$M_I\in End^{(i_1-i_0)}(\FF)\circ End^{(i_2-i_1)}(\FF)\circ \cdots \circ  End^{(i_k-i_{k-1})}(\FF) \subset End^{(i_k-i_0)}(\FF) \ .$$
Therefore,  Lemma \ref{le:partial decomposition} guarantees that
$\partial^{(i,j)}\in End^{(j-i)}(\FF)$
for all  $1\le i\le j\le \ell$.
%Indeed, taking in Lemma $\partial_i=\partial_{m_i}$ and $k=0$, we see that $\partial^{(i,j)}\in End^{(j-i)}(\FF)$.
Proposition \ref{pr:degree k noncommutative operator} is proved.
\end{proof}

We are ready now to finish the proof of Theorem \ref{th:difference derivative
implications}.

It is enough to apply Proposition \ref{pr:degree k noncommutative
operator} to the identity (\ref{eq:from delta to d}) from Lemma \ref{le: Delta(j-i)}.
\end{proof}

Finally, we are able to finish the proof of Theorem \ref{th:sl2 in sln}.

Note  that, according to Theorem \ref{th:difference derivative
implications}, the inclusions \eqref{eq:difference derivative}
imply the inclusions \eqref{eq:particular stabilization1}. Lemma
\ref{le:mij} and Proposition \ref{pr:mij and mij*} show that
\eqref{eq:particular stabilization1} and \eqref{eq:particular
stabilization2} are equivalent (one can see it by replacing $f_i$
by $f_{n-i-1}$, i.e., passing from from $m_i$ to $m_{n-i}^{-1}$ for all $i$).
Now the proof follows from Proposition \ref{pr:dij homogeneous}.

Theorem \ref{th:sl2 in sln} is proved.
\end{proof}

We will finish the section with a natural (yet conjectural) generalization of Theorem \ref{th:sl2 in sln}.

\begin{conjecture} Let $\gg=sl_2(\kk)$ and $A=A_n$ be as in Theorem \ref{th:sl2 in sln}, and $\FF$ be an object of ${\bf Alg}_1$.
 Then a matrix $g\in GL_n(\FF)$ belongs to
$G_{\gg,A}(\FF)$ if and only if
\begin{equation}
\label{eq:N-group characterization}
g \cdot \gg \cdot g^{-1} \subset (\gg,A)(\FF) \ .
\end{equation}

\end{conjecture}

\begin{remark} More generally, we would expect that for any perfect pair $(\gg,A)$ an element
$g\in (\FF\cdot A)^\times$ belongs to $G_{\gg,A}(\FF)$ if and only if \eqref{eq:N-group characterization} holds.
\end{remark}

\end{document}